\documentclass{article}
\usepackage{graphicx} 
\usepackage{amsbsy, amsthm, amsfonts, amsmath, amssymb, mathbbol, mathrsfs, enumerate, epsfig, graphicx, graphics, subfigure, tikz, MnSymbol, rotating, bm, multirow, color, xcolor, bbding, booktabs, threeparttable, appendix, natbib, tabularx}

\newtheorem{assumption}{Assumption}
\newtheorem{theorem}{Theorem}

\newtheorem{example}{Example}
\newtheorem{lemma}{Lemma}
\newtheorem{proposition}{Proposition}

\def\pr{\textup{pr}}
\def\PN{\textup{PN}}
\def\PC{\textup{PC}}

\def\T{\textup{T}}

\begin{document}

\parindent=22pt
\normalsize \baselineskip=18pt
\renewcommand{\baselinestretch}{1.2}
\renewcommand{\arraystretch}{1.5}
\renewcommand{\thefootnote}{\fnsymbol{footnote}}

\title{Identifying and bounding the probability of necessity for causes of effects with ordinal outcomes}

\author{
Chao Zhang$^1$, Zhi Geng$^1$, Wei Li$^2$, and Peng Ding$^3$ \\
$^1$ School of Mathematics and Statistics, Beijing Technology and Business University; \\
$^2$ Center for Applied Statistics and School of Statistics, Renmin University of China; \\
$^3$ Department of Statistics, University of California, Berkeley
}
\date{}
\maketitle

\begin{abstract}
Although the existing causal inference literature focuses on the forward-looking perspective by estimating effects of causes, the backward-looking perspective can provide insights into causes of effects. In backward-looking causal inference, the probability of necessity measures the probability that a certain event is caused by the treatment given the observed treatment and outcome. 
Most existing results focus on binary outcomes.
Motivated by applications with ordinal outcomes, we propose a general definition of the probability of necessity.
However, identifying the probability of necessity is challenging because it involves the joint distribution of the potential outcomes.
We propose a novel assumption of monotonic incremental treatment effect to identify the probability of necessity with ordinal outcomes.
We also discuss the testable implications of this key identification assumption. When it fails, we derive explicit formulas of the sharp large-sample bounds on the probability of necessity.
\end{abstract}


\section{Introduction to backward-looking causal inference and probability of necessity} \label{section: introduction}

Causal inference can be roughly classified into two categories: the forward-looking perspective that evaluates the effects of causes and the backward-looking perspective that evaluates the causes of effects. These are distinct statistics questions {\citep{pearl2015, dawid2022effects}}. The causal inference literature focuses on the forward-looking perspective by estimating, for example, the average causal effect, with randomized controlled trials as the gold standard for that task. By contrast, the backward-looking perspective asks different questions on how much we can attribute the effect to the treatment given the observed values of the treatment and outcome. With binary outcomes, \citet{pearl1999} proposed the notion of the probability of necessity to quantify the cause of effect, and discussed the identifiability and bounds of the probability of necessity with experimental data.  \citet{tianpearl2000} extended \citet{pearl1999} by utilizing both experimental data and observational data. \citet{kuroki2011statistical} used covariates to improve the bounds. For statistical inference of the probability of necessity with binary outcomes, \citet{cai2005variance} discussed frequentists' large-sample variances and  \citet{dawid2016statistical} applied the Bayesian method.

We focus on the backward-looking perspective for causal inference. In particular, we focus on ordinal outcomes, which are common in empirical research. We first propose the general definition of the probability of necessity with ordinal outcomes and illustrate its meaning with examples. With ordinal outcomes, the recent causal inference literature has made some progress from the forward-looking perspective \citep{ju2010, LuDingDasgupta2018, lu2020, gabriel2023}, but the corresponding discussion from the backward-looking perspective is missing. Our paper contributes to this gap.

However, evaluating causes of effects is more challenging than evaluating effects of causes, because the former inevitably involves the joint distribution of the potential outcomes \citep{pearl2015, dawid2022effects}. Due to the fundamental problem of causal inference, we can never jointly observe the potential outcomes which challenges the identification of their joint distribution. To identify the probability of necessity with binary outcomes, the conventional assumption is monotonicity that the treatment does not decrease the outcome at the individual level \citep{pearl1999, tianpearl2000}. Unfortunately, the monotonicity assumption does not suffice to identify the probability of necessity with ordinal outcomes. To overcome the challenge for identification, we propose a novel assumption of monotonic incremental treatment effect that the treatment does not decrease the outcome and increases the outcome by at most one level. We also derive testable implications of this key identification assumption. Moreover, when it is violated, we derive sharp large-sample bounds on the probability of necessity with ordinal outcomes.

\section{Potential outcomes and backward-looking causal inference with ordinal outcomes} \label{section: definition}

\subsection{Defining the probability of necessity with ordinal outcomes}
\label{section::p-o-n-ordinal}

Let $Z$ denote a binary treatment variable with 1 for treatment and 0 for control. Let $Y$ denote an ordinal outcome with $J$ levels labeled as $0, \ldots J-1$, where $0$ and $J-1$ represent the lowest and highest categories, respectively. 
Let $Y_1$ and $Y_0$ denote the potential outcomes under treatment and control, respectively. 
For a binary outcome, \cite{pearl1999} defined the probability of necessity as $\PN =\pr(Y_0=0 \mid  Z=1, Y=1)$ to measure how necessary the treatment $Z=1$ is for the occurrence of the outcome $Y=1$. Below we extend the PN to the case of an ordinal outcome. Let $\omega_0$ denote an event that depends only on the control potential outcome $Y_0$, e.g., $\omega_0= \{ Y_0 < y\}$ for some $y \in \{ 1, \ldots, J-1 \}$. To measure how necessary the treatment $Z=1$ is for an occurred outcome $Y=y$ changed from the counterfactual event $\omega_0$ without treatment, we define the probability of necessity as 
\[
\PN(\omega_0, y) = \pr(\omega_0 \mid Z = 1, Y = y).
\]
The $\PN(\omega_0, y)$ depends on the joint probability of $Y_1$ and $Y_0$ conditional on $Z = 1$. Let
\[
    q_{k\ell \mid z} = \pr(Y_1 = k, Y_0 = \ell \mid Z = z), \quad (z=0,1)
\]
summarized by the probability matrix $Q_z = (q_{k\ell \mid z})_{k, \ell = 0,\ldots, J-1}$. Then $\PN(\omega_0, y)$ is a function of $Q_1$ as
\begin{equation} \label{eq: defPN}
    \PN(\omega_0,y)
   =\frac{  \pr(\omega_0, Y_1 = y \mid Z=1) }{ \pr(Y = y \mid Z = 1) }
    = \frac{\sum_{\ell=0}^{J-1} c_\ell q_{y\ell \mid 1}}{\pr(Y = y \mid Z = 1)},
\end{equation}
where the $c_\ell$'s are binary and are uniquely determined by $\omega_0$. Therefore, we can identify $\PN(\omega_0, y)$ if we can identify $Q_1$. We illustrate the $\PN(\omega_0, y)$ through three examples below.

\begin{example} \label{ex: omega1}
    The Likert scale is a psychometric response scale primarily used in questionnaires to measure participants' preferences or degree of agreement with statements. 
    It typically employs a five-point scale ranging from strongly disagree on one end to strongly agree on the other, with neither agree nor disagree in the middle.
    Each level is assigned a numeric value or coding, usually starting at 0 and incremented by 1 for each level.
    Consider a participant's attitude score towards a certain question as outcome $Y$.
    Then $Y = 2$ means the neutrality.
    If $\omega_0 = \{Y_0 \neq 2\}$, then $\PN(\omega_0, 2) = \pr(Y_0 \neq 2 \mid Z = 1, Y = 2)$ measures the probability that a specific treatment $Z = 1$ is a necessary cause for individuals with $(Z = 1, Y = 2)$ to be neutral toward the question.
    From \eqref{eq: defPN}, we have $c_2 = 0$ and $c_0 = c_1 = c_3 = c_4 = 1$ for $\PN(Y_0 \neq 2, 2)$.
    We can identify $\PN(Y_0 \neq 2, 2)$ if we can identify $q_{2\ell \mid 1}$ $(\ell = 0, \ldots, 4)$.
    In general, with $\omega_0^1 = \{Y_0 \neq y\}$, define 
    \[
        \PN(\omega_0^1, y) = \pr(Y_0 \neq y \mid Z = 1, Y = y)
        = \frac{\sum_{\ell \neq y} q_{y\ell \mid 1}}{\pr(Y = y \mid Z = 1)} 
    \]
as the probability that $Z=1$ is the necessary cause of the observed $Y = y$.
\end{example}

\begin{example} \label{ex: omega2}
    \citet{oenema2001web} conducted a randomized controlled trial to assess whether web-based nutrition education changed personal attitude.
    Let $Y = 0, 1, 2$ when personal attitude is negative, neutral, and positive, respectively.
    If $\omega_0 = \{Y_0 = 0\}$, then $\PN(\omega_0, 2) = \pr(Y_0 = 0 \mid Z = 1, Y = 2)$ measures the probability that web-based nutrition education is the cause of individuals changing from negative attitude to the observed positive attitude.
    From \eqref{eq: defPN}, we have $c_0 = 1$ and $c_1 = c_2 = 0$ for $\PN(Y_0 = 0, 2)$.
    We can identify $\PN(Y_0 = 0, 2)$ if we can identify $q_{20 \mid 1}$.
    In general, with $\omega_0^2 = \{Y_0 = y'\}$, define
    \[
        \PN(\omega_0^2, y)
        = \pr(Y_0 = y' \mid Z = 1, Y = y)
        = \frac{ q_{yy' \mid 1}}{\pr(Y = y \mid Z = 1)}
    \]
      as the probability that $Z=1$ is the necessary cause of the observed outcome $Y = y$ changing from level $y'$.
\end{example}

\begin{example} \label{ex: omega3}
    \citet{lalonde1986} assessed the effect of a job training program on earnings.
    Let $Z$ be the binary indicator for receiving the job training program, and $Y \in \{0, 1, 2\}$ denote income levels with 0 for no income, 1 for low income, and 2 for high income.
    If $\omega_0 = \{Y_0 < 2\}$, then $\PN(\omega_0, 2) = \pr(Y_0 < 2 \mid Z = 1, Y = 2)$ measures the probability that job training is a necessary cause for individuals with observed high income to have increased their income from a lower level.
    From \eqref{eq: defPN}, we have $c_0 = c_1 = 1$ and $c_2 = 0$ for $\PN(Y_0<2, 2)$.
    In general, with $\omega_0^3 = \{Y_0 < y\}$, we define
    \[
        \PN(\omega_0^3, y)
        = \pr(Y_0 < y \mid Z = 1, Y = y)
        = \frac{\sum_{\ell < y} q_{y\ell \mid 1}}{\pr(Y = y \mid Z = 1)}
    \]
       as the probability that $Z=1$ is the necessary cause for the outcome to increase to the observed $Y = y$. 
\end{example}

\subsection{Comparing the probability of necessity and causal effects conditional on different evidence}

We focus on the probability of necessity $\PN(\omega_0, y)$ for backward-looking causal inference. An alternative strategy is to use the posterior causal effects conditional on the observed treatment and outcome \citep{lu2023, Li2023}, which provide finer information about treatment effect heterogeneity compared with the average effects. As it turns out, posterior causal effects are often special cases of $\PN(\omega_0, y)$.

We first consider the case with a binary $Y$. For forward-looking causal inference, the average causal effect, $\tau = E(Y_1 - Y_0)$, and the average causal effect on the treated units, $\tau_1 = E(Y_1 - Y_0 \mid Z=1)$, are commonly-used for evaluating effects of causes. For backward-looking causal inference, \citet{lu2023} and \citet{Li2023} proposed the posterior causal effect $E(Y_1 - Y_0 \mid Z=1, Y=1)$ conditional on the observed treatment $Z=1$ and outcome $Y=1$. We can verify that it equals $\pr(Y_0 = 0 \mid  Z=1, Y=1)$, which is a special case of $\PN(\omega_0^2, y)$ in Example \ref{ex: omega2}. The posterior causal effect $E(Y_1 - Y_0 \mid Z=1, Y=1)$ can have different sign compared with $\tau$ and $\tau_1$ in the presence of treatment effect heterogeneity across the subpopulations of $Y_1=1$ and $Y_1=0$. Therefore,   $\PN(\omega_0, y)$ and $E(Y_1 - Y_0 \mid Z=1, Y=1)$ can provide information beyond $\tau$ and $\tau_1$.

We then consider the case with an ordinal outcome $Y$. For forward-looking causal inference, \cite{ju2010} defined the distributional causal effect $\pr(Y_1 \geq y) - \pr(Y_0 \geq y)$, and \cite{LuDingDasgupta2018} defined $ \pr(Y_1 > Y_0)$ to measure the probability that the treatment is beneficial. For backward-looking causal inference, we can define the parallel posterior causal effects as $\pr(Y_1 \geq y  \mid  Z=1, Y=y) - \pr(Y_0 \geq y  \mid  Z=1, Y=y)$ and $\pr(Y_1 > Y_0  \mid Z=1, Y=y)$ conditional on the observed treatment $Z=1$ and outcome $Y=y$. We can verify that they are identical and both reduce to $\PN(Y_0<y, y) = \pr( Y_0 < y  \mid  Z=1, Y=y)$, which is a special case of $\PN(\omega_0^3, y)$ in Example \ref{ex: omega3}. The $\pr(Y_1 \geq y \mid  Z=1, Y=y) - \pr(Y_0 \geq y \mid  Z=1, Y=y)$ and $\pr(Y_1 > Y_0  \mid Z=1, Y=y)$ can  have different signs compared with $\pr(Y_1 \geq y) - \pr(Y_0 \geq y)$ and $ \pr(Y_1 > Y_0)$, respectively, when the treatment effects vary across subpopulations with $Y_1=y$ for different $y$'s.

\section{Identification and Partial Identification of the Probability of Necessity} \label{section: result}

\subsection{Identification of the probability of necessity for ordinal outcomes} \label{subsection: identify}

The definition of $\PN(\omega_0, y)$ involves the joint distribution $\pr(Y_1, Y_0 \mid Z=1)$. The observed data allows for the identification of the marginal distribution $\pr(Y_1\mid Z=1) = \pr(Y\mid Z=1)$ but not the marginal distribution $\pr(Y_0\mid Z=1)$ without further assumptions. To simplify the presentation, we assume $\pr(Y_0\mid Z=1)$ is identified and focus on identification of $\PN(\omega_0, y)$ with known marginal distributions of $\pr(Y_z\mid Z=1)$ for $z=1,0$.

\begin{assumption}
[identifiability of the counterfactual distribution]
\label{assume::marginal}
$\pr(Y_0\mid Z=1)$ is identifiable. 
\end{assumption}

Assumption \ref{assume::marginal} is strong but standard, which requires identifiability of the counterfactual distribution $\pr(Y_0\mid Z=1)$. The literature focuses on two sufficient conditions for Assumption \ref{assume::marginal}. First, under the unconfoundedness assumption that $Y_0$ is independent of $Z$ given observed covariates $X$, we can identify 
$$
\pr(Y_0=y \mid Z=1) 
= \sum_x \pr(Y_0=y\mid Z=1, x) \pr(x\mid Z=1)
= \sum_x \pr(Y=y\mid Z=0,x) \pr(x\mid Z=1),
$$
with ``$X=x$'' simplified as ``$x$.'' 
Second, \cite{tianpearl2000} focused on the setting with external experimental data which allows for identification of $\pr(Y_0)$ by $\pr_{\textup{E}}(Y\mid Z=0)$ due to randomization of $Z$, with the subscript ``E'' signifying the distribution under the experiment. Therefore, we can identify 
\begin{equation}\label{eq::experimental}
\pr(Y_0=y \mid Z=1) 
= \frac{ \pr(Y_0=y) - \pr(Z = 0, Y_0=y)  }{  \pr(Z = 1)  }
= \frac{\pr_{\textup{E}}(Y=y \mid Z = 0) - \pr(Z=0, Y=y)}{\pr(Z = 1)}. 
\end{equation}

Even under Assumption \ref{assume::marginal}, $\PN(\omega, y)$ is still not identifiable without additional assumptions. 
A commonly used assumption is monotonicity as stated in Assumption \ref{assume::monotonicity} below.

\begin{assumption}[monotonicity] \label{assume::monotonicity}
$Y_0 \leq Y_1$.
\end{assumption}

Assumption \ref{assume::monotonicity} requires that at the individual level, the treatment $Z$ does not have a negative effect on the outcome $Y$. For example, in the case of job training and income, Assumption \ref{assume::monotonicity} states that individuals who received job training would not have earned a higher income if they had not received the job training. Assumption \ref{assume::monotonicity} implies that $\pr(Y>y \mid Z=1,x) \geq \pr(Y>y \mid Z=0,x)$ for all $y$, which is a testable assumption based on the observed data. However, the observed data cannot validate Assumption \ref{assume::monotonicity} because it may not hold even if $\pr(Y>y \mid Z=1,x) \geq \pr(Y>y \mid Z=0,x)$ for all $y$.

Assumptions \ref{assume::marginal} and \ref{assume::monotonicity} are enough for identifying $\PN =\pr(Y_0=0  \mid  Z=1, Y=1)$ with binary outcomes \citep{pearl1999} but not enough with ordinal outcomes. Recall that the observed data can identify the marginal distribution $\pr(Y_1\mid Z=1)$ and Assumption \ref{assume::marginal} can identify the marginal distribution $\pr(Y_0\mid Z=1)$, which impose the following $2(J-1)$ linearly independent constraints on the $q_{k\ell \mid 1}$'s
\begin{equation}\label{eq: constraints PN}
\begin{aligned}
    &\sum_{\ell=0}^{J-1} q_{k\ell \mid 1} = \pr(Y_1 = k \mid Z = 1), \quad
     \sum_{k=0}^{J-1} q_{k\ell \mid 1} = \pr(Y_0 = \ell \mid Z = 1), \\
    &\sum_{k=0}^{J-1}\sum_{\ell=0}^{J-1} q_{k\ell \mid 1} = 1, \quad
      q_{k\ell \mid z} \geq 0 \quad (z = 0, 1;\ k, \ell = 0, \ldots, J-1).
\end{aligned}
\end{equation}
When $Y$ is binary with $J=2$, the joint probabilities are identifiable under the monotonicity assumption. With $J > 2$, the number of the unknown parameters in the joint probability of $Y_1$ and $Y_0$ in $Q_1$ equals $J^2-1$, which is larger than the number of equations in \eqref{eq: constraints PN}.

To identify $\PN(\omega_0, y)$ with ordinal outcome, below we make an assumption of monotonic incremental treatment effect which is stronger than the monotonicity assumption.

\begin{assumption}[Monotonic incremental treatment effect] \label{assume::incremental}
The treatment $Z$ does not decrease the level of $Y$ and increases $Y$ by at most one level, that is, $0 \leq Y_1 - Y_0 \leq 1$.
\end{assumption}

When $J = 2$, Assumption \ref{assume::incremental} reduces to the classic monotonicity assumption of $Y_1 \geq Y_0$.
Assumption \ref{assume::incremental} states that the treatment is harmless and can increase the outcome by at most one level for all units.
For example, it holds in Example \ref{ex: omega3} if the job training does not have a negative effect on income and limits the increase in income to no more than one level for all units.
Assumption \ref{assume::incremental} implies $q_{k,\ell \mid 1} = 0$ for $k < \ell$ or $k > \ell + 1$, which reduces the dimension of unknown parameters in $Q_1$ to $2J-2$. Then the $q_{k\ell \mid 1}$'s are identifiable by the constraints in \eqref{eq: constraints PN}. 
Lemma \ref{lemma::identification-Q} below summarizes the result with the definition of
\begin{equation} \label{eq: delta}
    \begin{aligned}
        \delta_{k \mid 1} = \sum_{j=0}^{k-1} \left\{ \pr(Y_0 = j \mid Z = 1) - \pr(Y_1 = j \mid Z = 1)\right\}. 
    \end{aligned}
    \end{equation}

\begin{lemma} \label{lemma::identification-Q}
    Under Assumptions \ref{assume::marginal} and \ref{assume::incremental}, $Q_1$ is identified by
    \begin{align*}
	  q_{00 \mid 1} = \pr(Y = 0 \mid Z = 1), \quad
        q_{k,k-1 \mid 1} = \delta_{k \mid 1},\quad
        q_{kk \mid 1} = \pr(Y = k \mid Z=1) - \delta_{k \mid 1}, \quad
        q_{k'\ell' \mid 1} = 0,
     \end{align*}
    for $k = 0, \ldots, J - 1$ and $k' < \ell'$ or $k' > \ell + 1$.
\end{lemma}

Lemma \ref{lemma::identification-Q} then ensures the identifiability of $\PN(\omega_0, y)$ by the relationship in \eqref{eq: defPN}.

\begin{theorem} \label{theorem::identification}
    Under Assumptions \ref{assume::marginal} and \ref{assume::incremental},
    $\PN(\omega_0, y)$ is identifiable by
    \[
	\PN(\omega_0, y) = c_{y} + \frac{ (c_{y-1} - c_y)\delta_{y \mid 1}}{\pr(Y = y \mid Z = 1)},
    \]
    where $c_{y-1}$ and $c_y$ are uniquely determined by $\omega_0$ in \eqref{eq: defPN}, and $\delta_{y \mid 1}$ is defined in \eqref{eq: delta}.
\end{theorem}

For sanity checks, we apply Theorem \ref{theorem::identification} to the case of binary $Y$ with $J = 2$. For example, let $\omega_0 = \{Y_0 = 0\}$, we have $c_0 = 1$ and $c_1 = 0$ for $\PN(Y_0 = 0, 1) = \pr(Y_0 = 0 \mid Z = 1, Y = 1)$.
Under Assumptions \ref{assume::marginal} and \ref{assume::incremental}, 
if $\pr(Y_0=0 \mid Z=1)$ is identified by external experimental data as in \eqref{eq::experimental}, then Theorem \ref{theorem::identification} recovers \citet{tianpearl2000}'s identification result on the probability of necessity for binary outcomes: 
\begin{equation*}
    \PN(Y_0 = 0, 1) 
    =  \frac{ \pr(Y_0 = 0 \mid Z = 1) - \pr(Y_1 = 0\mid Z=1)}{\pr(Y = 1 \mid  Z = 1 )} 
    = \frac{ \pr_{\textup{E}}(Y = 0 \mid Z = 0) - \pr(Y = 0)}{\pr(Z = 1, Y = 1)} .
\end{equation*}

\subsection{Sharp bounds on the probability of necessity for ordinal outcomes} \label{subsection: bounds}

Assumption \ref{assume::incremental} is key to the identification result in Theorem \ref{theorem::identification}. It has testable implications based on Lemma \ref{lemma::identification-Q} because the $q_{k\ell \mid 1}$'s must satisfy the Fr\'echet bounds \citep{ruschendorf1991frechet}. We present the result in Proposition \ref{proposition::falsify} below.

\begin{proposition} \label{proposition::falsify}
 Under Assumption \ref{assume::marginal}, Assumption \ref{assume::incremental} implies 
    \begin{align*}
		\max \left(\begin{array}{c}
			0 \\
			\pr(Y_1 = k \mid Z = 1) + \pr(Y_0 = k - 1 \mid Z = 1) - 1
		\end{array}\right)
        \leq &\delta_{k \mid 1} \leq
        \min \left(\begin{array}{c}
			\pr(Y_1 = k \mid Z = 1) \\
			\pr(Y_0 = k-1 \mid Z = 1)
		\end{array}\right), 
    \end{align*}
     for $k = 1, \ldots, J-1$ and $z = 0, 1$, where $\delta_{k \mid 1}$ is defined in \eqref{eq: delta}.
\end{proposition}

If the inequality in Proposition \ref{proposition::falsify} fails, then Assumption \ref{assume::incremental} is falsified. 
However, Assumption \ref{assume::incremental} cannot be validated by data. Even if the inequality in Proposition \ref{proposition::falsify} holds, Assumption \ref{assume::incremental} can still fail.
When Assumption \ref{assume::incremental} fails, the joint probability matrix $Q_1$ is not identifiable. Nevertheless, we can still derive sharp bounds on $\PN(\omega_0, y)$ as shown in Theorem \ref{theorem::bounds} below.

\begin{theorem} \label{theorem::bounds}
    Under Assumption \ref{assume::marginal}, the sharp bounds on $\PN(\omega_0, y)$ are
    \begin{align*}
		\max \left(\begin{array}{c}
			0,
			\frac{\pr(Y = y \mid Z = 1) - \sum_{\ell = 0}^{J-1} (1 - c_\ell) \pr(Y_0 = \ell \mid Z = 1)}{\pr(Y = y \mid Z = 1)}
		\end{array}\right)
        \leq \PN(\omega_0, y) \leq
        \min \left(\begin{array}{c}
			1,
			\frac{\sum_{\ell = 0}^{J-1} c_\ell \pr(Y_0 = \ell \mid Z = 1)}{\pr(Y = y \mid Z = 1)}\\
		\end{array}\right),
    \end{align*}
    where $c_l$'s are determined by $\omega_0$ in \eqref{eq: defPN}.
\end{theorem}

For sanity checks, we apply Theorem \ref{theorem::bounds} to the case of binary $Y$ with $J = 2$. Let $\omega_0 = \{Y_0 = 0\}$, we have $c_0 = 1$ and $c_1 = 0$ for $\PN(Y_0 = 0, 1) = \pr(Y_0 = 0 \mid Z = 1, Y = 1)$.
Under Assumption \ref{assume::marginal}, if $\pr(Y_0 \mid Z=1)$ is identified by external experimental data as in \eqref{eq::experimental}, then Theorem \ref{theorem::bounds} recovers \citet{tianpearl2000}'s bounds on the probability of necessity for binary outcomes: 
   \begin{align*}
	\max \left(\begin{array}{c}
		  0,
		  \frac{\pr(Y = 1) - \pr_{\textup{E}}(Y = 1 \mid Z = 0)}{\pr(Z = 1, Y = 1)}\\
	\end{array}\right)
	\leq \PN(Y_0 = 0, 1) \leq
	    \min \left(\begin{array}{c}
		  1,
		  \frac{ \pr_{\textup{E}}(Y = 0 \mid Z = 0) - \pr(Z = 0, Y = 0)}{\pr(Z = 1, Y = 1)}
	\end{array}\right). 
    \end{align*}

Theorem \ref{theorem::bounds} drops Assumption \ref{assume::incremental} entirely. If we are willing to maintain the monotonicity in Assumption \ref{assume::monotonicity}, then we can derive narrower sharp bounds on $\pr(Y_0 \neq y \mid Z = 1, Y = y)$ and $\pr(Y_0 = y' \mid Z = 1, Y = y)$.

\begin{theorem}
\label{theorem::bounds-mono}
Under Assumptions \ref{assume::marginal} and \ref{assume::monotonicity},  
the sharp bounds on $\pr(Y_0 \neq y \mid Z = 1, Y = y)$ are    
$$
\max \left(\begin{array}{c}
			0,
			\frac{\pr(Y = y \mid Z = 1) - \pr(Y_0 = y \mid Z = 1)}{\pr(Y = y \mid Z = 1)}
		\end{array}\right)
        \leq \pr(Y_0 \neq y \mid Z = 1, Y = y) \leq
        \min \left(\begin{array}{c}
			1,
			\frac{\delta_{y \mid 1}}{\pr(Y = y \mid Z = 1)}
		\end{array}\right),
$$
and the sharp bounds on $\pr(Y_0 = y' \mid Z = 1, Y = y)$ are 
\begin{align*}
		&\max \left(\begin{array}{c}
			0,
			\frac{\pr(Y = y \mid Z = 1) + \sum_{k=0}^{y'-1} \pr(Y = k \mid Z = 1) - \sum_{\ell = 0}^y \pr(Y_0 = \ell \mid Z = 1) + \pr(Y_0 = y' \mid Z = 1)}{\pr(Y = y \mid Z = 1)}
		\end{array}\right)\\
        &\leq \pr(Y_0 = y' \mid Z = 1, Y = y) \leq
        \min_{y' < k \leq y} \left(\begin{array}{c}
			1,
            \frac{\pr(Y_0 = y' \mid Z = 1)}{\pr(Y = y \mid Z = 1)},
			\frac{\delta_{k \mid 1}}{\pr(Y = y \mid Z = 1)}
		\end{array}\right),
  \end{align*}
  where $\delta_{k \mid 1}$ is defined in \eqref{eq: delta}. 
\end{theorem}

Theorem \ref{theorem::bounds} provides the sharp bounds on general $\PN(\omega_0, y)$ under Assumption \ref{assume::incremental}, whereas Theorem \ref{theorem::bounds-mono} provides the sharp bounds on two special cases of $\PN(\omega_0, y)$ under the weaker Assumption \ref{assume::monotonicity}. Under monotonicity in Assumption \ref{assume::monotonicity}, the explicit formulas for the sharp bounds on general $\PN(\omega_0, y)$ do not have simple elegant formulas. Nevertheless, we can still calculate the bounds numerically by solving linear programming problems. We relegate the computational details to the Supplementary Materials.

\section{Illustration} \label{section: application}

We illustrate the theory with the Lalonde data reviewed in Example \ref{ex: omega3} with both the experimental source that contains 445 samples and the observational source that contains 16,177 samples. The data were used to study causal effects of job training programs on income. 
We discretize the outcome into three categories with 0 for zero income, 1 for medium income below 2,650 dollars, and 2 for high income above 2,650 dollars.
We apply one-to-one matching based on the observed covariates in the observational data to ensure better overlap of covariates under treatment and control groups.  Table \ref{application: data} summarizes the data. For simplicity of presentation, we focus on point estimates below and omit the results based on confidence intervals.

\begin{table}[t]
    \centering
    \caption{The Lalonde data with both the experimental source and the matched observational source}
    {\begin{tabular}{cccccccc}
    &        &\multicolumn{3}{c}{experimental data}   &\multicolumn{3}{c}{observational data}\\
    &        &   $Y=0$   &   $Y=1$   &   $Y=2$        &   $Y=0$   &   $Y=1$   &   $Y=2$      \\
    &$Z=1$   &     45    &     32    &    108         &     90    &     64    &    216       \\
    &$Z=0$   &     92    &     33    &    135         &     115   &     50    &    205       \\
    \end{tabular}
    }
    \label{application: data}
\end{table}

As benchmark analyses, we first estimate the forward-looking causal quantities based on the experimental data. The estimates for the distributional causal effects are $\pr(Y_1 \geq 1) - \pr(Y_0 \geq 1) = 0.12$ and $\pr(Y_1 \geq 2) - \pr(Y_0 \geq 2) = 0.07$, which indicate that the job training increases the proportion of both $Y=1$ and $Y=2$. We also obtain the estimates for the bounds $0.20 \leq \pr(Y_1>Y_0)-\pr(Y_1 < Y_0) \leq 0.29$ and $0.58 \leq \pr(Y_1 > Y_0) \leq 1$, which also indicate that the job training is beneficial. Overall, forward-looking causal inference suggests positive effects of the job training on income.

\begin{table}[t]
   \centering
   \caption{Identification of $\PN(\omega_0, y)$ by Theorem \ref{theorem::identification}, and sharp bounds on $\PN(\omega_0,y)$ by Theorems \ref{theorem::bounds} and \ref{theorem::bounds-mono}, respectively 
    }
    {\resizebox{\textwidth}{!}{
    \begin{tabular}{ccccccccc}
         	               & $\PN(\omega_0, y)$                       &$\PN(Y_0 \neq y, y)$   &$\PN(Y_0 =0, y)$    &$\PN(Y_0=1, y)$    &$\PN(Y_0=2, y)$    &$\PN(Y_0 < y, y)$ \\
    \multirow{3}{*}  {$y=2$}  & Theorem \ref{theorem::identification}    & 0.17                  & 0.00               &  0.17             &  0.83             &  0.17  	         \\
                              & Theorem \ref{theorem::bounds}            &[0.17, 0.88]           &[0.00, 0.68]        &[0.00, 0.20]       &[0.12, 0.83]       &[0.17, 0.88]      \\
                              & Theorem \ref{theorem::bounds-mono}       &[0.17, 0.17]           &[0.00, 0,17]        &[0.00, 0.17]       &[0.83, 0,83]       &[0.17, 0.17]      \\
    \multirow{3}{*}  {$y=1$}  & Theorem \ref{theorem::identification}    & 0.89                  & 0.89               &  0.11             &  0.00             & 0.89     	     \\
                              & Theorem \ref{theorem::bounds}            &[0.31, 1.00]           &[0.00, 1.00]        &[0.00, 0.69]       &[0.00, 1.00]       &[0.00, 1.00]      \\
                              & Theorem \ref{theorem::bounds-mono}       &[0.31, 0.89]           &[0.31, 0.89]        &[0.11, 0.69]       &[0.00, 0.00]       &[0.31, 0.89]      \\
    \end{tabular} }
    } 
    \label{application: result}
\end{table}

For backward-looking causal inference, we report the estimates and bounds on various $\PN(\omega_0,y)$ in Table \ref{application: result}. We first discuss the estimates by Theorem \ref{theorem::identification} under Assumptions \ref{assume::marginal} and  \ref{assume::incremental}. Assumption \ref{assume::marginal} holds because we use the experimental data to identify $\pr(Y_0\mid Z=1)$ by \eqref{eq::experimental}. Assumption \ref{assume::incremental} requires that the job training does not harm income and increases income by at most one level. We highlight two estimates in  Table \ref{application: result}. The $(1,4)$th cell contains the estimate $\PN(Y_0=2,2)=\pr(Y_0=2 \mid Z=1, Y=2)=0.83$. It shows that individuals with high income after receiving the job training have 83\% probability of having high income even without receiving the job training. The $(4,2)$th cell contains the estimate $\PN(Y_0=0,1)=\pr(Y_0=0 \mid Z=1, Y=1)=0.89$. It shows that individuals with medium income after receiving the job training have 89\% probability of having zero income without receiving the job training. Overall, the backward-looking causal inference suggests that it is more likely that the treatment effect occurs for individuals having medium income, not high income, after receiving the job training.

Dropping Assumption \ref{assume::incremental} entirely, $\PN(\omega_0,y)$ is not identifiable and the bounds by Theorem \ref{theorem::bounds} can be wide as shown in Table \ref{application: result}. 
Nevertheless, if we are willing to maintain the monotonicity under Assumption \ref{assume::monotonicity}, many bounds by Theorem \ref{theorem::bounds-mono} shrink and even reduce to the point estimates. 
The bounds in the $(3,4)$th cell collapse to the same point estimate $\PN(Y_0=2,2)=\pr(Y_0=2 \mid Z=1, Y=2)=0.83$. The bounds in the $(6,2)$th cell are $0.31 \leq \PN(Y_0=0,1)=\pr(Y_0=0 \mid Z=1, Y=1) \leq 0.89$, with the upper bound $0.89$ being the point estimate and the lower bound $0.31$ still providing evidence that individuals with medium income after receiving the job training have at least 31\% probability of having zero income without receiving the job training.
The qualitative results remain the same as above. Again, compared with the forward-looking causal inference,  the backward-looking causal inference provides finer information about heterogeneous treatment effects conditional on the observed information of the treatment and outcome.

\section{Discussion of the Probability of Causation} \label{section: dissussion}

To measure the likelihood that $Z = 1$ is a cause of an effect $Y = 1$ with a binary $Y$, \citet{dawidfaigmanfienberg2014} defined the probability of causation as $\PC = \pr\left(Y_0 = 0 \mid Y_{1} = 1\right)$ and discussed its identifiability and bounds; see also \citet{robins1989probability} and \citet{greenland1999relation}. The PC only considers whether the treatment is a cause of changing the potential outcomes, without considering whether the target individuals received the treatment or not. The PC is equivalent to $\PN =\pr(Y_0=0  \mid  Z=1, Y=1)$ when there is no confounding between $Z$ and $Y$, that is, $Z$ is independent of $(Y_1, Y_0)$. However, in observational studies, the PC differs from the PN if $Z $ is not independent of $(Y_1, Y_0)$. With an ordinal outcome, we can define the general probability of causation as $\PC(\omega_0, y) = \pr(\omega_0 \mid Y_1 = y)$ in parallel with $\PN(\omega_0, y) = \pr(\omega_0 \mid Z=1,  Y = y)$. We can derive the identification and bounding results on $\PC(\omega_0, y)$ in parallel with Theorems \ref{theorem::identification} -- \ref{theorem::bounds-mono}. We provide details in the Supplementary Material.



\section*{Supplementary material} \label{SM}

\appendix
\setcounter{theorem}{0} 
\setcounter{lemma}{0} 
\setcounter{assumption}{0} 
\setcounter{proposition}{0} 
\setcounter{example}{0} 
\setcounter{table}{0} 
\renewcommand{\thetheorem}{S\arabic{theorem}}
\renewcommand{\thelemma}{S\arabic{lemma}}
\renewcommand{\theassumption}{S\arabic{assumption}}
\numberwithin{equation}{section}
\renewcommand{\theproposition}{S\arabic{proposition}} 
\renewcommand{\theexample}{S\arabic{example}}
\renewcommand{\thetable}{S\arabic{table}}

The supplementary material contains the following sections. 

Section \ref{section: lemma} proves Lemma \ref{lemma::identification-Q} and illustrates it with a numerical example. 

Section \ref{section: theorem identify} proves Theorem \ref{theorem::identification}.

Section \ref{section: theorem bound} proves Theorem \ref{theorem::bounds} in three steps. 
Section \ref{subsection: bound} states a lemma and proves the bounds. Section \ref{subsection: sharp lower} and Section \ref{subsection: sharp upper} prove the sharpness of lower and upper bounds, respectively. 

Section \ref{section: bounds-mono} provides the sharp bounds on $\PN(\omega_0,y)$ under monotonicity. 
Sections \ref{subsection: bounds-mono-proof}--\ref{subsection: bounds-mono-upper} prove Theorem \ref{theorem::bounds-mono}. 
Section \ref{subsection: linear programming} provides a linear programming method to derive sharp bounds on general $\PN(\omega_0,y)$ under monotonicity. 

Section \ref{section: proposition} proves Proposition \ref{proposition::falsify}. 

Section \ref{section: discussion} provides the identification results and sharp bounds on the probability of causation with ordinal outcomes.

\section{Proof of lemma \ref{lemma::identification-Q}} \label{section: lemma}

In this section, we prove Lemma \ref{lemma::identification-Q} and illustrate it with Example \ref{example::illustrate-lemma1} below.

\begin{proof}
    Under Assumption \ref{assume::marginal}, the marginal distribution of $Y_1$ and that of $Y_0$ conditional on $Z = 1$ are identifiable and provide $2(J-1)$ linearly independent constraints on $q_{k\ell \mid 1}$’s as shown in \eqref{eq: constraints PN}. 
    Assumption \ref{assume::incremental} further imposes the following constraints on $q_{k\ell \mid 1}$'s: 
    \begin{equation} \label{A.eq: assumption2}
        q_{k\ell \mid 1} = 0, 
    \end{equation} 
    for $k < \ell$ or $k > \ell+1$.
    Therefore, the dimension of unknown parameters in the joint probability matrix $Q_1$ reduces to $2J-2$. 
    By \eqref{eq: constraints PN} and \eqref{A.eq: assumption2}, the constraints on $q_{k\ell \mid 1}$'a are  
    \begin{equation}\label{A.eq: constraints Q_1}
    \begin{aligned}
	q_{00 \mid 1} &= \pr(Y_1 = 0 \mid Z = 1),\\
        q_{k,k-1 \mid 1} + q_{kk \mid 1} &= \pr(Y_1 = k \mid Z = 1), \\
	q_{k-1,k-1 \mid 1} + q_{k,k-1 \mid 1} &= \pr(Y_0 = k-1 \mid Z = 1), \\
        q_{k'\ell' \mid 1} &= 0,
    \end{aligned}
    \end{equation}
    for $k = 1, \ldots, J-1$ and $k' < \ell'$ or $k' > \ell' + 1$.

    Solving the system of linear equations in \eqref{A.eq: constraints Q_1}, we can identify all $q_{k\ell \mid 1}$'s in $Q_1$ as 
    \begin{align*}
	  q_{00 \mid 1} &= \pr(Y = 0 \mid Z = 1), \\
        q_{k,k-1 \mid 1} &= \sum_{j=0}^{k-1} \left\{ \pr(Y_0 = j \mid Z = 1) - \pr(Y = j \mid Z = 1)\right\},\\ 
         q_{kk \mid 1} &= \sum_{j=0}^{k} \pr(Y = j \mid Z = 1) - \sum_{j=0}^{k-1} \pr(Y_0 = j \mid Z = 1), \\
        q_{k'\ell' \mid 1} &= 0, 
     \end{align*}
    for $k = 1, \ldots, J-1$ and $k' < \ell'$ or $k' > \ell'+1$. 
\end{proof}

Then, we illustrate Lemma \ref{lemma::identification-Q} with Example \ref{example::illustrate-lemma1} below. 
\begin{example}\label{example::illustrate-lemma1}
    Let the probability matrix $Q_1 = (q_{k\ell \mid 1})_{k,\ell= 0,\ldots, J-1}$ summarize the joint distribution of the potential outcomes conditional on $Z=1$. 
    For a case with $J = 5$, the probability matrix is 
    \[
    Q_1 =\  \bordermatrix{
             ~        & q_{+0 \mid 1} & q_{+1 \mid 1} & q_{+2 \mid 1} & q_{+3 \mid 1} & q_{+4 \mid 1} \cr
        q_{0+ \mid 1} & q_{00 \mid 1} & q_{01 \mid 1} & q_{02 \mid 1} & q_{03 \mid 1} & q_{04 \mid 1} \cr
        q_{1+ \mid 1} & q_{10 \mid 1} & q_{11 \mid 1} & q_{12 \mid 1} & q_{13 \mid 1} & q_{14 \mid 1} \cr
        q_{2+ \mid 1} & q_{20 \mid 1} & q_{21 \mid 1} & q_{22 \mid 1} & q_{23 \mid 1} & q_{24 \mid 1} \cr
        q_{3+ \mid 1} & q_{30 \mid 1} & q_{31 \mid 1} & q_{32 \mid 1} & q_{33 \mid 1} & q_{34 \mid 1} \cr
        q_{4+ \mid 1} & q_{40 \mid 1} & q_{41 \mid 1} & q_{42 \mid 1} & q_{43 \mid 1} & q_{44 \mid 1} 
    }, 
    \]
    where $q_{k+ \mid 1} = \sum_{\ell=0}^{4} q_{k\ell \mid 1} = \pr(Y_1 = k \mid Z = 1)$ for $k = 0,1,2,3,4$ and $q_{+\ell \mid 1} = \sum_{k=0}^{4} q_{k\ell \mid 1} = \pr(Y_0 = \ell \mid Z = 1)$ for $\ell = 0,1,2,3,4$ are the marginal probabilities. 
    The dimension of the unknown free parameters is $5\times 5 - 1 = 24$. 
    Suppose Assumption \ref{assume::marginal} holds. We can identify $q_{k+ \mid 1}$ and $q_{+\ell \mid 1}$, which provide 8 linearly independent constraints on $Q_1$. Under Assumption \ref{assume::incremental}, we have 
    \[
    Q_1 =\  \bordermatrix{
             ~        & q_{+0 \mid 1} & q_{+1 \mid 1} & q_{+2 \mid 1} & q_{+3 \mid 1} & q_{+4 \mid 1} \cr
        q_{0+ \mid 1} & q_{00 \mid 1} &       0       &       0       &       0       &       0       \cr
        q_{1+ \mid 1} & q_{10 \mid 1} & q_{11 \mid 1} &       0       &       0       &       0       \cr
        q_{2+ \mid 1} &       0       & q_{21 \mid 1} & q_{22 \mid 1} &       0       &       0       \cr
        q_{3+ \mid 1} &       0       &       0       & q_{32 \mid 1} & q_{33 \mid 1} &       0       \cr
        q_{4+ \mid 1} &       0       &       0       &       0       & q_{43 \mid 1} & q_{44 \mid 1} 
    }. 
    \]
    The number of the unknown free parameters is reduced to 8, which indicates that the number of linearly independent constraints provided by $q_{k+ \mid 1}$ and $q_{+\ell \mid 1}$ is equal to the number of unknown parameters. In fact, there is a unique solution for $q_{k\ell \mid 1}$, so we can identify $Q_1$. 
\end{example}

\section{Proof of Theorem \ref{theorem::identification}} \label{section: theorem identify}

When Assumptions \ref{assume::marginal} and \ref{assume::incremental} hold, we can identify $Q_1$ by Lemma \ref{lemma::identification-Q}. 
Because $\PN(\omega_0,y)$ is a linear function of the $q_{k\ell \mid 1}$'s, we have 
\[
    \begin{aligned}
        \PN(\omega_0,y) 
       =&\frac{\sum_{\ell=0}^{J-1} c_\ell q_{y\ell \mid 1}}{\pr(Y = y \mid Z = 1)} \\
       =& \frac{c_{y-1} q_{y,y-1 \mid 1} + c_y q_{yy \mid 1}}{\pr(Y = y \mid Z = 1)} \\
       =& c_y + \frac{\sum_{j=0}^{k-1} \left\{ \pr(Y_0 = j \mid Z = 1) - \pr(Y = j \mid Z = 1)\right\}}{\pr(Z = 1,Y = y)}, 
    \end{aligned}
\]
where $c_{y-1}$ and $c_y$ are uniquely determined by $\omega_0$ in \eqref{eq: defPN}.

\section{Proof of Theorem \ref{theorem::bounds}} \label{section: theorem bound}

\subsection{Proving the bounds} \label{subsection: bound}

In this subsection, we show 
\begin{align*}
    &\PN_{\textup{L}}(\omega_0,y)
    = \max \left(\begin{array}{c}
	   0 \\
	   \frac{\pr(Y = y \mid Z = 1) - \sum_{\ell=0}^{J-1} (1 - c_\ell) \pr(Y_0 = \ell \mid Z = 1)}{\pr(Y = y \mid Z = 1)}
	 \end{array}\right),\\
    &\PN_{\textup{U}}(\omega_0,y) 
    = \min \left(\begin{array}{c}
	     1 \\
	     \frac{\sum_{\ell=0}^{J-1} c_\ell \pr(Y_0 = \ell \mid Z = 1)}{\pr(Y = y \mid Z = 1)}
      \end{array}\right)
\end{align*}
are the lower and upper bounds on $\PN(\omega_0,y)$, respectively.

We first recall a lemma.

\begin{lemma}[Fr\'echet bounds \citep{ruschendorf1991frechet}] \label{lemma: frechet}
    For any three events A, B and C, we have 
    \begin{align*}
	\max \left\{0, \pr(A \mid C) + \pr(B \mid C) - 1 \right\} 
	\leq \pr(A, B \mid C) 
	\leq \min \left\{ \pr(A \mid C), \pr(B \mid C) \right\}.
    \end{align*}
\end{lemma}

Then, we prove $\PN_{\textup{L}}(\omega_0,y)$ and $\PN_{\textup{U}}(\omega_0,y)$ are the lower and upper bounds on $\PN(\omega_0,y)$, respectively.

\begin{proof} 
    For any given $\omega_0$ and observed evidence $(Z = 1, Y = y)$, we have 
    \begin{equation}\label{C.eq: defPN}
    \begin{aligned}
         \PN(\omega_0,y) 
        =\frac{\pr(Y_1 = y, \omega_0 \mid Z = 1)}{\pr(Y = y \mid Z = 1)}. 
    \end{aligned}
    \end{equation}
    In \eqref{C.eq: defPN}, $\pr(Y = y \mid Z = 1)$ is directly identified from the observed data. 
    We only derive the sharp bounds on $\pr(Y_1 = y, \omega_0 \mid Z = 1)$ under the constraints in \eqref{eq: constraints PN}.

    By the Fr\'echet bounds in Lemma \ref{lemma: frechet}, we have
    \begin{align}
	\pr(Y_1 = y, \omega_0 \mid Z = 1) 
	&\geq \max\left\{0, \pr(Y_1 = y \mid Z = 1) + \pr(\omega_0 \mid Z = 1) - 1\right\} ,  \label{C.eq: Lfrechet}\\
        \pr(Y_1 = y,\omega_0 \mid Z = 1) 
	&\leq \min\left\{\pr(Y_1 = y \mid Z = 1), \pr(\omega_0 \mid Z = 1)\right\}. \label{C.eq: Ufrechet}
    \end{align}
    Then, we have 
    \begin{align}
	\PN(\omega_0,y) 
        \geq &\frac{\max\{0, \pr(Y_1 = y \mid Z = 1) + \pr(\omega_0 \mid Z = 1) - 1\}}{\pr(Y = y \mid Z = 1)} , \label{C.eq: LPN}\\
	  \PN(\omega_0, y) 
	\leq &\frac{\min\{\pr(Y_1 = y \mid Z = 1), \pr(\omega_0 \mid Z = 1)\}}{\pr(Y = y \mid Z = 1)}. \label{C.eq: UPN}
    \end{align}
    From the definition of $\omega_0$, we have $\pr(\omega_0 \mid Z = 1) = \sum_{\ell=0}^{J-1} c_\ell \pr(Y_0 = \ell \mid Z = 1)$. 
    Therefore, by \eqref{C.eq: LPN} and \eqref{C.eq: UPN}, we have  
    \begin{align*}
        \PN(\omega_0, y)
	& \geq \max \left(\begin{array}{c}
	   0 \\
	   \frac{\pr(Y = y \mid Z = 1) - \sum_{\ell=0}^{J-1} (1 - c_\ell) \pr(Y_0 = \ell \mid Z = 1)}{\pr(Y = y \mid Z = 1)}
	 \end{array}\right)
            = \PN_{\textup{L}}(\omega_0, y), \\
	\PN(\omega_0,y)
	& \leq \min \left(\begin{array}{c}
	     1 \\
	     \frac{\sum_{\ell=0}^{J-1} c_\ell \pr(Y_0 = \ell \mid Z = 1)}{\pr(Y = y \mid Z = 1)}
      \end{array}\right)
            = \PN_{\textup{U}}(\omega_0, y).
    \end{align*} 
\end{proof}

\subsection{Proving the sharpness of the lower bound} \label{subsection: sharp lower}

In this subsection, we prove the sharpness of the lower bound by showing the existence of a matrix $Q_1$ satisfying \eqref{eq: constraints PN} such that $\PN(\omega_0,y)$ attains $\PN_{\textup{L}}(\omega_0, y)$. 
We first state a lemma used in the construction of $Q_1$.

\begin{lemma} \label{lemma: existence}
    For any given $n$-dimensional vector $q_1=(q_{1+}, \ldots, q_{n+})^\T$ and $m$-dimensional vector $q_0=(q_{+1}, \ldots, q_{+m})^\T$, if $q_{k+} \geq 0$ for $k=1, \ldots, n$, $q_{+\ell} \geq 0$ for $\ell=1, \ldots, m$, and $\sum_{k=1}^n q_{k+} = \sum_{\ell=1}^m q_{+\ell}$, then under the constraints of row sums $q_1$ and column sums $q_0$, there exists at least one non-negative matrix $Q_{n \times m}$.
\end{lemma}

We prove Lemma \ref{lemma: existence} by demonstrating the existence of $Q$ under the constraints of $q_1$ and $q_0$.

\begin{proof}[of Lemma \ref{lemma: existence}]
    Let 
    \[
        q_{k\ell} = \frac{q_{k+} q_{+\ell}}{S}, \quad 
        (k=1,\ldots, n; \ell = 1, \ldots, m), 
    \]
    denote the element in the $k$th row and $\ell$th column of $Q$, where $S = \sum_{k=1}^n q_{k+} = \sum_{\ell=1}^m q_{+\ell}$. 
    Next, we show that such $Q$ satisfies the constraints of $q_1$ and $q_0$. 
    
    First, $q_{k\ell} \geq 0$ because $q_{k+} \geq 0$ for $k=1, \ldots, n$ and $q_{+\ell} \geq 0$ for $\ell = 1, \ldots, m$. 
    Then, the sum of elements in the $k$th row of matrix $Q$ is 
    \[
        \sum_{\ell=1}^m q_{k\ell}
        = \sum_{\ell=1}^m \frac{q_{k+} q_{+\ell}}{S} 
        = \frac{q_{k+}}{S} \sum_{\ell=1}^m q_{+\ell} 
        = q_{k+}, 
        \quad (k = 1, \ldots, n). 
    \]
    The sum of elements in the $\ell$th column of matrix $Q$ is 
    \[
        \sum_{k=1}^n q_{k\ell}
        = \sum_{k=1}^n \frac{q_{k+} q_{+\ell}}{S} 
        = \frac{q_{+\ell}}{S} \sum_{k=1}^n q_{k+} 
        = q_{+\ell}, 
        \quad (\ell = 1, \ldots, m). 
    \]
    That is, $Q=(q_{k\ell})_{k=1,\ldots,n;\ell=1,\ldots,m}$ is a non-negative matrix that satisfies the constraints of row sums $q_1$ and column sums $q_0$, when $q_{k+} \geq 0$ for $k=1, \ldots, n$, $q_{+\ell} \geq 0$ for $\ell = 1, \ldots, m$ and $\sum_{k=1}^n q_{k+} = \sum_{\ell=1}^m q_{+\ell}$. 
\end{proof}

Next, we prove the existence of a matrix $Q_1$ satisfying \eqref{eq: constraints PN} such that $\PN(\omega_0,y)$ attains the lower bound.

 \begin{proof}
    We can represent the marginal distribution of potential outcomes conditional on $Z = 1$ by the row sums and column sums of $Q_1$ as 
    \[
        q_{k+ \mid 1} = \sum_{\ell=0}^{J-1} q_{k\ell \mid 1} = \pr(Y_1 = k \mid Z = 1), \quad 
	  q_{+\ell \mid 1} = \sum_{k=0}^{J-1} q_{k\ell \mid 1}=\pr(Y_0 = \ell \mid Z = 1), \quad
        \left(k, \ell = 0, \ldots, J-1 \right). 
    \]
    Let $q_{1 \mid 1} = \left(q_{0+ \mid 1}, \ldots, q_{J-1,+ \mid 1}\right)^{\T}$ and 
    $q_{0 \mid 1} = \left(q_{+0 \mid 1}, \ldots, q_{+,J-1 \mid 1}\right)^{\T}$. 
    Under Assumption \ref{assume::marginal}, we can identify $q_{1 \mid 1}$ and $q_{0 \mid 1}$. 
    By \eqref{C.eq: Lfrechet}, $Q_1$ makes the lower bound attainable implies that if $\pr(Y_1 = y \mid Z = 1) + \pr(\omega_0 \mid Z = 1) - 1 \leq 0$, we have 
    \begin{align} \label{C.eq: Lbound1}
        \pr(Y_1 = y, \omega_0 \mid Z = 1) = 0. 
    \end{align}
    Otherwise, 
    \begin{align} \label{C.eq: Lbound2}
        \pr(Y_1 = y, \omega_0 \mid Z = 1) = \pr(Y_1 = y \mid Z = 1) + \pr(\omega_0 \mid Z = 1) - 1. 
    \end{align}

    \noindent {\bf Step 1}: 
    If $\pr(Y_1 = y \mid Z = 1) + \pr(\omega_0 \mid Z = 1) - 1 \leq 0$, $Q_1$ is subject to \eqref{eq: constraints PN} and \eqref{C.eq: Lbound1}. 
    By \eqref{C.eq: Lbound1}, we have 
	\begin{equation*}\label{C.eq: proof Lbound1}
	\begin{aligned}
		  \pr(Y_1 = y, \omega_0 \mid Z = 1)
            = \sum_{\ell=0}^{J-1} c_\ell \pr(Y_1 = y,Y_0 = \ell \mid Z = 1) 
            = \sum_{\ell=0}^{J-1} c_\ell q_{y\ell \mid 1}
            = 0. 
	\end{aligned}
	\end{equation*}
    This implies that all elements in the $y$th row corresponding to columns in $Q_1$ where $c_\ell = 1$ should be zero. 
    Then, we construct the elements of the $y$th row in $Q_1$ as 
    \[
	q_{y0 \mid 1} = \min\left\{(1 - c_0)q_{y+ \mid 1}, q_{+0 \mid 1}\right\}, \quad 
	q_{y\ell \mid 1} = \min\left\{(1 - c_\ell)\left(q_{y+ \mid 1} - \sum_{j=0}^{\ell-1} q_{yj \mid 1}\right), q_{+\ell \mid 1}\right\}, 
    \]
    for $\ell = 1, \ldots, J-1$. 
    Let $Q^*_1 = (q^*_{k\ell})_{k=0, \ldots, J-2; \ell=0, \ldots, J-1}$ denote the submatrix of $Q_1$ composed of all rows except the $y$th row. 
    Then, the row sums and column sums of the submatrix $Q_1^*$ are 
    \begin{align*}
	q^*_{1} &= \left(q^*_{0+}, \ldots, q^*_{J-2,+ }\right)^\T
	=\left(q_{0+ \mid 1}, \ldots,q_{y-1,+ \mid 1}, q_{y+1,+ \mid 1}, \ldots, q_{J-1,+ \mid 1}\right)^\T, \\
	q^*_{0} &= \left(q^*_{+0}, \ldots, q^*_{+,J-1}\right)^\T
	= \left(q_{+0 \mid 1} - q_{y0 \mid 1}, \ldots, q_{+,J-1 \mid 1} - q_{y,J-1 \mid 1}\right)^\T, 
    \end{align*}
    which satisfy $q^*_{k+} \geq 0$ for $k=0,\ldots,J-2$ and $q^*_{+\ell} \geq 0$ for $\ell=0,\ldots,J-1$. 
    Because $\pr(Y_1 = y \mid Z = 1) + \pr(\omega_0 \mid Z = 1) - 1 \leq 0$, by \eqref{eq: constraints PN}, we have 
    \[
    \sum_{k=0}^{J-2} q^*_{k+} = \sum_{\ell=0}^{J-1} q^*_{+\ell}= \sum_{k=0}^{J-2}\sum_{\ell=0}^{J-1} q^*_{k\ell}= 1 - q_{y+ \mid 1}. 
    \]
    Apart from $q^*_{1}$ and $q^*_{0}$, there are no additional constraints on $Q_1^*$. 
    By Lemma \ref{lemma: existence}, at least one matrix $Q^*_1$ satisfies the constraints of row sums $q^*_1$ and column sums $q^*_0$ which implies the existence of $Q_1$.  
    That is, there exists a probability matrix $Q_1$ subject to constraints from observed data such that $\PN(\omega_0,y)$ attains $\PN_{\text{L}}(\omega_0, y) = 0$ when $\pr(Y_1 = y \mid Z = 1) + \pr(\omega_0 \mid Z = 1) - 1 \leq 0$.

    \noindent {\bf Step 2:} If $\pr(Y_1 = y \mid Z = 1) + \pr(\omega_0 \mid Z = 1) -1 > 0$, $Q_1$ is subject to \eqref{eq: constraints PN} and \eqref{C.eq: Lbound2}. 
    By \eqref{C.eq: Lbound2}, we have 
    \begin{equation*}\label{C.eq: proof Lbound2}
	\pr(Y_1 \neq y, \bar{\omega}_0 \mid Z = 1) 
	= \sum_{\ell=0}^{J-1} (1 - c_\ell) \pr(Y_1 \neq y, Y_0 = \ell \mid Z = 1) 
	= \sum_{k\neq y}\sum_{\ell=0}^{J-1} (1 - c_\ell) q_{k\ell \mid 1} 
        = 0,
    \end{equation*}
    where $\bar{\omega}_0$ is the complement set of $\omega_0$.
    This implies 
    \[
	q_{y\ell \mid 1} = q_{+\ell \mid 1}, 
    \]
    for $c_\ell = 0$ where $\ell = 0, \ldots, J-1$, and
    \[
        q_{k\ell \mid 1} = 0, 
    \]
    for $c_\ell = 0$ where $k \neq y$ and $k, \ell = 0, \ldots, J-1$. 
    That is all elements in columns where $c_\ell=0$ are determined.

    Let $Q'_1$ denote a submatrix formed by the undetermined elements in $Q_1$. 
    Then, the row sums and column sums of the submatrix $Q'_1$ are 
    \begin{align*}
        q'_1 &= \left(q'_{0+},\ldots,q'_{J-1,+}\right)^\T
             = \left(q_{0+ \mid 1}, \ldots, q_{y-1,+ \mid 1}, 
                q_{y+ \mid 1} - \sum_{\ell=0}^{J-1} (1 - c_\ell) q_{+\ell \mid 1}, 
                q_{y+1,+ \mid 1}, \ldots, q_{J-1,+ \mid 1}\right)^\T, \\
	q'_0 &= \left(q'_{+0}, \ldots, q'_{+,\sum_{\ell=0}^{J-1} c_\ell}\right)^\T,
    \end{align*}
    where $q'_{0}$ is formed by arranging the elements of 
    $(c_0q_{+0 \mid 1}, \ldots, c_{J-1}q_{+,J-1 \mid 1})^\T$ for $c_{\ell} = 1$, $\ell = 0, \ldots, J-1$ in sequential order, 
    and $q'_{k+} \geq 0$ for $k = 0, \ldots, J-1$ and $q'_{+\ell} \geq 0$ for $\ell = 0, \ldots, \sum_{\ell=0}^{J-1} c_\ell$.

    Because $\pr(Y_1 = y \mid Z = 1) + \pr(\omega_0 \mid Z = 1) - 1 > 0$, by \eqref{eq: constraints PN}, we have 
    \[
	\sum_{k=0}^{J-1} q^\prime_{k+}
	= \sum_{\ell=0}^{\sum_{\ell=0}^{J-1} c_\ell} q^\prime_{+\ell}
	= \sum_{k=0}^{J-1}\sum_{\ell=0}^{\sum_{\ell=0}^{J-1} c_\ell} q^\prime_{k\ell}
	= 1- \sum_{\ell=0}^{J-1} (1-c_\ell) q_{+\ell \mid 1}.
    \]
    Apart from $q'_{1}$ and $q'_{0}$, there are no additional constraints on $Q_1'$. 
    By Lemma \ref{lemma: existence}, at least one matrix $Q'_1$ satisfies the constraints of row sums $q'_1$ and column sums $q'_0$ which implies the existence of $Q_1$. 
    That is, there exists a probability matrix $Q_1$ subject to constraints from observed data such that $\PN(\omega_0,y)$ attains $\PN_{\text{L}}(\omega_0, y) = \pr(Y_1 = y \mid Z = 1) + \pr(\omega_0 \mid Z = 1) -1$ when $\pr(Y_1 = y \mid Z = 1) + \pr(\omega_0 \mid Z = 1) -1 > 0$. 
\end{proof}

\subsection{Proving the sharpness of the upper bound} \label{subsection: sharp upper}

In this subsection, we prove the sharpness of the upper bound by proving the existence of a matrix $Q_1$ satisfying \eqref{eq: constraints PN} such that $\PN(\omega_0,y)$ attains $\PN_{\textup{U}}(\omega_0, y)$.

\begin{proof}
    In \eqref{C.eq: Ufrechet}, $Q_1$ makes the upper bound attainable implies that if $\pr(Y_1 = y \mid Z = 1) \geq \pr(\omega_0 \mid Z = 1)$, we have 
    \begin{align} \label{C.eq: Ubound1}
        \pr(Y_1 = y, \omega_0 \mid Z = 1) = \pr(\omega_0 \mid Z = 1). 
    \end{align}
    Otherwise,  
    \begin{align} \label{C.eq: Ubound2}
        \pr(Y_1 = y, \omega_0 \mid Z = 1) = \pr(Y_1 = y \mid Z = 1). 
    \end{align}

\noindent {\bf Step 1:} If $\pr(\omega_0 \mid Z = 1) \leq \pr(Y_1 = y \mid Z = 1)$, $Q_1$ is subject to \eqref{eq: constraints PN} and \eqref{C.eq: Ubound1}. 
    By \eqref{C.eq: Ubound1}, we have
    \begin{equation*} \label{C.eq: proof Ubound1}
	    \pr(Y_1 \neq y, \omega_0 \mid Z = 1) 
         = \sum_{\ell=0}^{J-1} c_\ell \pr(Y_1 \neq y,Y_0 = \ell \mid Z = 1)
         = \sum_{k \neq y}\sum_{\ell=0}^{J-1} c_\ell q_{k\ell \mid 1}
         = 0. 
    \end{equation*}
    This implies 
    \[
	q_{y\ell \mid 1} = q_{+\ell \mid 1}, 
    \]
    for $c_\ell = 1$ where $\ell = 0, \ldots, J-1$, and
    \[
        q_{k\ell \mid 1} = 0,
    \]
    for $c_\ell = 1$ where $k \neq y$ and $k, \ell = 0, \ldots, J-1$. 
    That is all elements in columns where $c_\ell=1$ are determined.

    Let $Q_1^*$ denote a submatrix formed by the undetermined elements in $Q_1$. 
    Then, the row sums and column sums of the submatrix $Q^*_1$ are 
    \begin{align*}
	q^*_1 &= \left(q^*_{0+}, \ldots, q^*_{J-1,+}\right)^\T
         = \left(q_{0+ \mid 1}, \ldots, q_{y-1,+ \mid 1}, 
          q_{y+ \mid 1}-\sum_{\ell=0}^{J-1} c_\ell q_{+\ell \mid 1}, 
          q_{y+1,+ \mid 1}, \ldots, q_{J-1,+ \mid 1}\right)^\T, \\
	q^*_0 &= \left(q^*_{+0}, \ldots, q^*_{+,J-1-\sum_{\ell=0}^{J-1} c_\ell}\right)^\T,
    \end{align*}
    where $q_{0}^*$ is formed by arranging the elements of $((1-c_0)q_{+0 \mid 1}, \ldots, (1-c_{J-1})q_{+,J-1 \mid 1})^\T$ for $c_{\ell} = 0$, $\ell = 0, \ldots, J-1$ in sequential order, and $q^*_{k+} \geq 0$ for $k = 0, \ldots, J-1$ and $q^*_{+\ell} \geq 0$ for $\ell = 0, \ldots, J - 1 - \sum_{\ell=0}^{J-1} c_\ell$.

    Because $\pr(\omega_0 \mid Z = 1) \leq \pr(Y_1 = y \mid Z = 1)$, by \eqref{eq: constraints PN}, we have 
    \[
	\sum_{k=0}^{J-1} q^*_{k+}
	= \sum_{\ell=0}^{J-1-\sum_{\ell=0}^{J-1} c_\ell} q^*_{+\ell}
	= \sum_{k=0}^{J-1}\sum_{\ell=0}^{J-1-\sum_{\ell=0}^{J-1} c_\ell} q^*_{k\ell}
	= 1- \sum_{\ell=0}^{J-1} c_\ell q_{+\ell \mid 1}.
    \]
    Apart from $q^*_{1}$ and $q^*_{0}$, there are no additional constraints on $Q_1^*$. 
    By Lemma \ref{lemma: existence}, at least one matrix $Q^*_1$ satisfies the constraints of row sums $q^*_1$ and column sums $q^*_0$ which implies the existence of $Q_1$.  
    That is, there exists a probability matrix $Q_1$ subject to constraints from observed data such that $\PN(\omega_0,y)$ attains $\PN_{\textup{U}}(\omega_0,y) = \pr(\omega_0 \mid Z = 1)$ when $\pr(\omega_0 \mid Z = 1) \leq \pr(Y_1 = y \mid Z = 1)$.

\noindent {\bf Step 2:} If $\pr(\omega_0 \mid Z = 1) > \pr(Y_1 = y \mid Z = 1)$, $Q_1$ is subject to \eqref{eq: constraints PN} and \eqref{C.eq: Ubound2}. 
    By \eqref{C.eq: Ubound2}, we have
    \begin{equation*} \label{C.eq: proof Ubound2}
	\pr(Y_1 = y, \bar{\omega}_0 \mid Z = 1)
	  = \sum_{\ell=0}^{J-1} (1 - c_\ell)\pr(Y_1 = y,Y_0 = \ell \mid Z = 1) 
        = \sum_{\ell=0}^{J-1} (1 - c_\ell)q_{y\ell \mid 1}
        =0. 
    \end{equation*}
    This implies that all elements in the $y$th row corresponding to columns in $Q_1$ where $c_\ell = 0$ should be zero. 
    Then, we construct the elements of the $y$th row in matrix $Q_1$ as 
    \[
	q_{y0 \mid 1} = \min\left\{c_0q_{y+ \mid 1}, q_{+0 \mid 1}\right\}, \quad 
	q_{y\ell \mid 1} = \min\left\{c_\ell \left(q_{y+ \mid 1} - \sum_{j=0}^{\ell-1}q_{yj \mid 1}\right), q_{+\ell \mid 1}\right\}, \quad 
	(\ell = 1,\ldots, J-1).
    \]
    Let $Q_1' = (q'_{k\ell})_{k=0,\ldots, J-2; \ell=0, \ldots, J-1}$ denote the submatrix of $Q_1$ composed of all rows except the $y$th row. 
    Then, the row sums and column sums of the submatrix $Q_1'$ are 
    \begin{align*}
	q'_1 &= \left(q'_{0+}, \ldots, q'_{J-2,+}\right)^\T
	= \left(q_{0+ \mid 1}, \ldots, q_{y-1,+ \mid 1}, q_{y+1,+ \mid 1}, \ldots, q_{J-1,+ \mid 1}\right)^\T, \\
	q'_0 &= \left(q'_{+0}, \ldots, q'_{+,J-1}\right)^\T
	= \left(q_{+0 \mid 1} - q_{y0 \mid 1}, \ldots, q_{+,J-1 \mid 1} - q_{y,J-1 \mid 1}\right)^\T.
    \end{align*}
    which satisfy $q^*_{k+} \geq 0$ for $k=0,\ldots,J-2$ and $q^*_{+\ell} \geq 0$ for $\ell=0,\ldots,J-1$. 
    Because $\pr(\omega_0 \mid Z = 1) > \pr(Y_1 = y \mid Z = 1)$, by \eqref{eq: constraints PN}, we have 
    \[
	\sum_{k=0}^{J-2} q^\prime_{k+}
	= \sum_{\ell=0}^{J-1} q^\prime_{+\ell}
	= \sum_{k=0}^{J-2}\sum_{\ell=0}^{J-1} q^\prime_{k\ell}
	= 1- q_{y+ \mid 1}.
    \]

    Apart from $q'_{1}$ and $q'_{0}$, there are no additional constraints on $Q_1'$. 
    Similarly, by Lemma \ref{lemma: existence}, at least one matrix $Q'_1$ satisfies the constraints of row sums $q'_1$ and column sums $q'_0$ which implies the existence of $Q_1$. 
    That is, there exists a probability matrix $Q_1$ subject to constraints from observed data such that $\PN(\omega_0,y)$ attains 
    $\PN_{\textup{U}}(\omega_0,y) = \pr(Y_1 = y \mid Z = 1)$ when $\pr(\omega_0 \mid Z = 1) > \pr(Y_1 = y \mid Z = 1)$. 
\end{proof}

\section{Sharp bounds on $\PN(\omega_0,y)$ under monotonicity} \label{section: bounds-mono}

In this section, we first prove Theorem \ref{theorem::bounds-mono} in three subsections. Section \ref{subsection: bounds-mono-proof} proves the upper and lower bounds, Section \ref{subsection: bounds-mono-lower} proves that the lower bound is sharp, and Section \ref{subsection: bounds-mono-upper} proves that the upper bound is sharp. Then, Section \ref{subsection: linear programming} provides a linear programming method to derive sharp bounds on the general $\PN(\omega_0,y)$ under monotonicity.

\subsection{Proving the bounds in Theorem \ref{theorem::bounds-mono}} \label{subsection: bounds-mono-proof}

Recall the definition of $\delta_{k \mid 1}$ in \eqref{eq: delta}. 
In this subsection, we only prove that 
\begin{align*}
    \PN_{\text{U}}(Y_0 = y', y) &= 
    \min \left(\begin{array}{c}
			1,
            \frac{\pr(Y_0 = y' \mid Z = 1)}{\pr(Y = y \mid Z = 1)}, 
			\frac{\min_{y' < k \leq y} \delta_{k \mid 1}}{\pr(Y = y \mid Z = 1)}
		\end{array}\right), \\
    \PN_{\text{L}}(Y_0 = y', y) &= 
    \max \left(\begin{array}{c}
			0,
			\frac{\pr(Y = y \mid Z = 1) + \sum_{k=0}^{y'-1} \pr(Y = k \mid Z = 1) - \sum_{\ell = 0}^y \pr(Y_0 = \ell \mid Z = 1) + \pr(Y_0 = y' \mid Z = 1)}{\pr(Y = y \mid Z = 1)}
		\end{array}\right),
\end{align*}
are the sharp upper and lower bounds on $\PN(Y_0 = y', y)$ for $y' \leq y$, respectively, because $\pr(Y_0 \neq y \mid Z = 1, Y = y) = 1 - \pr(Y_0 = y \mid Z = 1, Y = y)$ under monotonicity.

Under Assumption \ref{assume::monotonicity}, the probability matrix $Q_1$ is a lower triangular matrix where $q_{k\ell \mid 1} = 0$ for $k<\ell$. 
The marginal distribution of potential outcomes conditional on $Z = 1$ provides the following constraints on $q_{k\ell \mid 1}$'s: 
\begin{equation} \label{E.eq: constraints PN}
    \begin{aligned}
        &q_{k+ \mid 1} = \sum_{\ell=0}^{k} q_{k\ell \mid 1}=\pr(Y_1=k \mid Z = 1), \quad \left(k=0, \ldots, J-1 \right) \\ 
        &q_{+\ell \mid 1} = \sum_{k=\ell}^{J-1} q_{k\ell \mid 1}=\pr(Y_0=\ell \mid Z = 1), \quad \left(\ell=0, \ldots, J-1 \right) \\
        &\sum_{k=\ell}^{J-1}\sum_{\ell=0}^{J-1} q_{k\ell \mid 1}=1,\quad  q_{k\ell \mid 1} \geq 0, \quad \left(k, \ell=0, \ldots, J-1 \right). 
    \end{aligned}
\end{equation}
Let $q_{1 \mid 1} = \left(q_{0+ \mid 1}, \ldots, q_{J-1,+ \mid 1}\right)^{\T}$ and 
    $q_{0 \mid 1} = \left(q_{+0 \mid 1}, \ldots, q_{+,J-1 \mid 1}\right)^{\T}$. 
Under Assumption \ref{assume::marginal}, we can identify $q_{1 \mid 1}$ and $q_{0 \mid 1}$. 
By \eqref{eq: defPN}, we have 
\[
    \pr(Y_0 = y' \mid Z = 1, Y = y)
    = \frac{\pr(Y_1 = y, Y_0 = y' \mid Z = 1)}{\pr(Y= y \mid Z = 1)}
    = \frac{q_{yy' \mid 1}}{\pr(Y=y \mid Z = 1)}
    = \frac{\sum_{\ell=0}^{J-1} c_\ell q_{y\ell \mid 1}}{\pr(Y = y \mid Z = 1)}, 
\]
where $c_{y'} = 1$ and $c_\ell = 0$ for $\ell \neq y'$. 
Since $\pr(Y=y \mid Z=1)$ is identified from observed data, we only derive sharp bounds on $\pr(Y_1=y, Y_0=y' \mid Z = 1)$ under \eqref{E.eq: constraints PN}, which is equal to proving that 
\begin{align*}
          & \max \left(0, \pr(Y = y \mid Z = 1) + \pr(Y_0 = y' \mid Z = 1) - \{\pr(Y_0 \leq y \mid Z = 1) - \pr(Y_1 < y' \mid Z = 1)\} \right) \\
     \leq & \pr(Y_1 = y, Y_0 = y' \mid Z = 1) \\
     \leq & \min \left(\pr(Y_1 = y \mid Z = 1), \pr(Y_0 = y' \mid Z = 1), \min_{y' \leq h < y} \{\pr(Y_0 \leq h \mid Z = 1) - \pr(Y_1 \leq h \mid Z = 1)\} \right),
\end{align*}
holds for $y' \leq y$ and the bounds are attainable under \eqref{E.eq: constraints PN}.

Next, we first provide the proof for $y'<y$ in two steps. 
In the first step, we prove that $\max (0, \pr(Y = y \mid Z = 1) + \pr(Y_0 = y' \mid Z = 1) - \{\pr(Y_0 \leq y \mid Z = 1) - \pr(Y_1 < y' \mid Z = 1)\} )$ is the lower bound on $\pr(Y_1 = y, Y_0 = y' \mid Z = 1)$ under monotonicity assumption. 
In the second step, we prove $\min_{y' \leq h < y} \{\pr(Y_1 = y \mid Z = 1), \pr(Y_0 = y' \mid Z = 1), \pr(Y_0 \leq h \mid Z = 1) - \pr(Y_1 \leq h \mid Z = 1) \}$ is the upper bound on $\pr(Y_1 = y, Y_0 = y' \mid Z = 1)$ under monotonicity assumption. 
Then, the proof when $y'=y$ is similarly provided.

\begin{proof}

\noindent {\bf Step 1:} Proving the lower bound for $y'<y$. 

    For the lower bound, $0 \leq \pr(Y_1 = y, Y_0 = y' \mid Z = 1)$ naturally holds. 
    We next prove $\pr(Y = y \mid Z = 1) + \pr(Y_0 = y' \mid Z = 1) - \{\pr(Y_0 \leq y \mid Z = 1) - \pr(Y_1 < y' \mid Z = 1)\} \leq \pr(Y_1 = y, Y_0 = y' \mid Z = 1)$. 
    \begin{align*}
         & \pr(Y = y \mid Z = 1) + \pr(Y_0 = y' \mid Z = 1) - \{\pr(Y_0 \leq y \mid Z = 1) - \pr(Y_1 < y' \mid Z = 1)\} \\
        =& \sum_{\ell=0}^{J-1} \pr(Y_1 = y, Y_0 =\ell \mid Z = 1) + \sum_{k=0}^{J-1} \pr(Y_1 = k, Y_0 = y' \mid Z = 1) \\
         &- \left\{ \sum_{\ell=0}^{y}\sum_{k=0}^{J-1} \pr(Y_1 = k, Y_0 = \ell \mid Z = 1) - \sum_{k=0}^{y'-1}\sum_{\ell=0}^{J-1} \pr(Y_1 = k, Y_0 = \ell \mid Z = 1) \right\} \\
        =& \sum_{\ell=0}^{y}\pr(Y_1 = y, Y_0 = \ell \mid Z = 1) + \sum_{k=y'}^{J-1} \pr(Y_1 = k, Y_0 = y' \mid Z = 1) - \sum_{k=y'}^{J-1}\sum_{\ell=0}^{y} \pr(Y_1 = k, Y_0 = \ell
         \mid Z = 1) \\
        =& \pr(Y_1 = y, Y_0 = y' \mid Z = 1) - \sum_{k=y',k \neq y}^{J-1}\sum_{\ell=0,\ell \neq y'}^{y} \pr(Y_1 = k, Y_0 = \ell \mid Z = 1) \\
        \leq& \pr(Y_1 = y, Y_0 = y' \mid Z = 1), 
    \end{align*}
    where the second equality holds because of Assumption \ref{assume::monotonicity}. 
    Therefore, we have 
    \begin{equation} \label{eq::D.lbound1}
    \begin{aligned}
        & \max \left(0, \pr(Y = y \mid Z = 1) + \pr(Y_0 = y' \mid Z = 1) - \{\pr(Y_0 \leq y \mid Z = 1) - \pr(Y_1 < y' \mid Z = 1)\} \right) \\
        \leq & \pr(Y_1 = y, Y_0 = y' \mid Z = 1). 
    \end{aligned}
    \end{equation}

\noindent {\bf Step 2:} Proving the upper bound for $y'<y$. 

    For the upper bound, we have $\pr(Y_1 = y, Y_0 = y' \mid Z = 1) \leq \min\left\{\pr(Y_1 = y \mid Z = 1), \pr(Y_0 = y' \mid Z = 1) \right\}$ from Lemma \ref{lemma: frechet}. 
    We next prove that $\pr(Y_1 = y, Y_0 = y' \mid Z = 1) \leq \min_{y' < h \leq y} \{\pr(Y_0 \leq h \mid Z = 1) - \pr(Y_1 \leq h \mid Z = 1)\}$. 
    \begin{align*}
         & \pr(Y_0 \leq h \mid Z = 1) - \pr(Y_1 \leq h \mid Z = 1) \\ 
        =& \sum_{k=0}^{J-1}\sum_{\ell=0}^{h} \pr(Y_1 = k, Y_0 = \ell \mid Z = 1) - \sum_{\ell=0}^{J-1}\sum_{k=0}^{h} \pr(Y_1 = k, Y_0 = \ell \mid Z = 1) \\
        =& \sum_{k=h}^{J-1}\sum_{\ell=0}^{h} \pr(Y_1 = k, Y_0 = \ell \mid Z = 1) \\
        =& \pr(Y_1 = y, Y_0 = y' \mid z = 1) +  \sum_{k=h,k\neq y}^{J-1}\sum_{\ell=0,\ell \neq y'}^{h} \pr(Y_1 = k, Y_0 = \ell \mid Z = 1) \\
        \geq& \pr(Y_1 = y, Y_0 = y' \mid z = 1), 
    \end{align*}
    where the second equality holds because of Assumption \ref{assume::monotonicity} and the last equality holds because of $y' \leq h < y$. 
    Therefore, we have 
    \begin{equation} \label{eq::D.ubound1}
    \begin{aligned}
      & \pr(Y_1 = y, Y_0 = y' \mid Z = 1) \\
     \leq & \min_{y' \leq h < y} \left\{\pr(Y_1 = y \mid Z = 1), \pr(Y_0 = y' \mid Z = 1), \pr(Y_0 \leq h \mid Z = 1) - \pr(Y_1 \leq h \mid Z = 1) \right\}. 
    \end{aligned}
    \end{equation}

\noindent {\bf Step 3:} Proving the lower bound for $y'=y$.

For the lower bound, $0 \leq \pr(Y_1 = y, Y_0 = y \mid Z = 1)$ naturally holds. 
    We next prove $\pr(Y = y \mid Z = 1) + \pr(Y_0 = y \mid Z = 1) - \{\pr(Y_0 \leq y \mid Z = 1) - \pr(Y_1 < y \mid Z = 1)\} \leq \pr(Y_1 = y, Y_0 = y \mid Z = 1)$. 
    \begin{align*}
         & \pr(Y = y \mid Z = 1) + \pr(Y_0 = y \mid Z = 1) - \{\pr(Y_0 \leq y \mid Z = 1) - \pr(Y_1 < y \mid Z = 1)\} \\
        =& \sum_{\ell=0}^{J-1} \pr(Y_1 = y, Y_0 =\ell \mid Z = 1) + \sum_{k=0}^{J-1} \pr(Y_1 = k, Y_0 = y \mid Z = 1) \\
         &- \left\{ \sum_{\ell=0}^{y}\sum_{k=0}^{J-1} \pr(Y_1 = k, Y_0 = \ell \mid Z = 1) - \sum_{k=0}^{y-1}\sum_{\ell=0}^{J-1} \pr(Y_1 = k, Y_0 = \ell \mid Z = 1) \right\} \\
        =& \sum_{\ell=0}^{y}\pr(Y_1 = y, Y_0 = \ell \mid Z = 1) + \sum_{k=y}^{J-1} \pr(Y_1 = k, Y_0 = y \mid Z = 1) - \sum_{k=y}^{J-1}\sum_{\ell=0}^{y} \pr(Y_1 = k, Y_0 = \ell
         \mid Z = 1) \\
        =& \pr(Y_1 = y, Y_0 = y \mid Z = 1) - \sum_{k=y+1}^{J-1}\sum_{\ell=0}^{y-1} \pr(Y_1 = k, Y_0 = \ell \mid Z = 1) \\
        \leq& \pr(Y_1 = y, Y_0 = y \mid Z = 1), 
    \end{align*}
    where the second equality holds because of Assumption \ref{assume::monotonicity}. 
    Therefore, we have 
    \begin{equation} \label{eq::D.lbound2}
    \begin{aligned}
        & \max \left(0, \pr(Y = y \mid Z = 1) + \pr(Y_0 = y \mid Z = 1) - \{\pr(Y_0 \leq y \mid Z = 1) - \pr(Y_1 < y \mid Z = 1)\} \right) \\
        \leq & \pr(Y_1 = y, Y_0 = y \mid Z = 1). 
    \end{aligned}
    \end{equation}

\noindent {\bf Step 4:} Proving the upper bound for $y'=y$.

    For the upper bound, we have $\pr(Y_1 = y, Y_0 = y \mid Z = 1) \leq \min\left\{\pr(Y_1 = y \mid Z = 1), \pr(Y_0 = y \mid Z = 1) \right\}$ from Lemma \ref{lemma: frechet}. 
    When $y' = y$, $\max_{y' \leq h <y} \{\pr(Y_0 \leq h \mid Z = 1) - \pr(Y_1 \leq h \mid Z = 1)\}$ is an empty set. 
    Therefore, 
    \begin{align} \label{eq::D.ubound2}
    \pr(Y_1 = y, Y_0 = y \mid Z = 1) \leq \min\{\pr(Y_1 = y \mid Z = 1), \pr(Y_0 = y \mid Z = 1)\} = \PN_{\textup{U}}(Y_0=y,y).    
    \end{align}    
\end{proof}

\subsection{Proving the sharpness of the lower bound in Theorem \ref{theorem::bounds-mono}} \label{subsection: bounds-mono-lower}

In this subsection, we prove the sharpness of the lower bounds in two steps.

\begin{proof}
\noindent{\bf Step 1:} Proving the sharpness of the lower bound for $y'<y$.

    By \eqref{eq::D.lbound1}, the attainability of the lower bound implies that there exists a joint probability matrix $Q_1$ satisfying \eqref{E.eq: constraints PN} such that equality holds. 
    This means that if $\pr(Y = y \mid Z = 1) + \pr(Y_0 = y' \mid Z = 1) - \{\pr(Y_0 \leq y \mid Z = 1) - \pr(Y_1 < y' \mid Z = 1) \} \leq 0$, we have 
    \begin{equation}\label{E.eq: Lbound1}
        \pr(Y_0=y', Y_1=y \mid Z=1) = 0.
    \end{equation}
    Otherwise, 
    \begin{align} \label{E.eq: Lbound2}
        &\pr(Y_0=y', Y_1=y \mid Z=1)\\ \nonumber
        &= \pr(Y = y \mid Z = 1) + \pr(Y_0 = y' \mid Z = 1) - \{\pr(Y_0 \leq y \mid Z = 1) - \pr(Y_1 < y' \mid Z = 1) \}. 
    \end{align}

    Then, we prove the sharpness of the lower bound by showing the existence of a matrix $Q_1$ satisfying \eqref{E.eq: constraints PN} such that $\Pr(Y_1=y, Y_0=y' \mid Z=1)$ attains the bound in two parts.

    \begin{itemize}
        \item [1.] If $\pr(Y = y \mid Z = 1) + \pr(Y_0 = y' \mid Z = 1) - \{\pr(Y_0 \leq y \mid Z = 1) - \pr(Y_1 < y' \mid Z = 1) \} \leq 0$, $Q_1$ is subject to \eqref{E.eq: constraints PN} and \eqref{E.eq: Lbound1}. 
        By \eqref{E.eq: Lbound1}, we have 
        \begin{equation*}\label{E.eq: proof Lbound1}
	\begin{aligned}
		  \pr(Y_1 = y,Y_0 = y' \mid Z = 1) 
            = q_{yy' \mid 1}
            = 0. 
	\end{aligned}
	\end{equation*}
        Then, we construct the elements of the $y$th row in $Q_1$ as
        \[
	       q_{y0 \mid 1} = \min\left\{(1 - c_0)q_{y+ \mid 1}, q_{+0 \mid 1}\right\}, \quad 
	       q_{y\ell \mid 1} = \min\left\{(1 - c_\ell)\left(q_{y+ \mid 1} - \sum_{j=0}^{\ell-1} q_{yj \mid 1}\right), q_{+\ell \mid 1}\right\} 
        \]
        for $\ell = 1, \ldots, y$. 
        Let $Q_1^* = (q^*_{k\ell})_{k=0, \ldots, J-2; \ell=0, \ldots, J-1}$ denote the submatrix of $Q_1$ composed of all rows except the $y$th row. 
        Then, the row sums and column sums of the submatrix $Q_1^*$ are 
        \begin{align*}
	   q^*_{1} &= \left(q^*_{0+}, \ldots, q^*_{J-2,+ }\right)^\T
    	=\left(q_{0+ \mid 1}, \ldots,q_{y-1,+ \mid 1}, q_{y+1,+ \mid 1}, \ldots, q_{J-1,+ \mid 1}\right)^\T, \\
	   q^*_{0} &= \left(q^*_{+0}, \ldots, q^*_{+,J-1}\right)^\T
    	= \left(q_{+0 \mid 1} - q_{y0 \mid 1}, \ldots, q_{+,J-1 \mid 1} - q_{y,J-1 \mid 1}\right)^\T,
        \end{align*}
        which satisfy $q^*_{k+} \geq 0$ for $k=0,\ldots,J-2$ and $q^*_{+\ell} \geq 0$ for $\ell=0,\ldots,J-1$. 
         Because $\pr(Y = y \mid Z = 1) + \pr(Y_0 = y' \mid Z = 1) - \{\pr(Y_0 \leq y \mid Z = 1) - \pr(Y_1 < y' \mid Z = 1) \} \leq 0$, we have 
         \[
         \sum_{k=0}^{J-2} q^*_{k+} = \sum_{\ell=0}^{J-1} q^*_{+\ell}= \sum_{k=0}^{J-2}\sum_{\ell=0}^{J-1} q^*_{k\ell \mid 1}= 1 - q_{y+ \mid 1}. 
         \]
        Apart from $q^*_{1}$ and $q^*_{0}$, there are no additional constraints on $Q_1^*$. 
        By Lemma \ref{lemma: existence}, at least one matrix $Q^*_1$ satisfies the constraints of row sums $q^*_1$ and column sums $q^*_0$ which implies the existence of $Q_1$.

        \item [2.] If $\pr(Y = y \mid Z = 1) + \pr(Y_0 = y' \mid Z = 1) - \{\pr(Y_0 \leq y \mid Z = 1) - \pr(Y_1 < y' \mid Z = 1) \} > 0$, $Q_1$ is subject to \eqref{E.eq: constraints PN} and \eqref{E.eq: Lbound2}. 
        In this part, we first construct the determinable elements in $Q_1$ based on \eqref{E.eq: constraints PN}, \eqref{E.eq: Lbound2} and the monotonicity assumption. 
        Then, we demonstrate that the undetermined elements of $Q_1$ satisfy \eqref{E.eq: constraints PN}, \eqref{E.eq: Lbound2}, and the monotonicity Assumption \ref{assume::monotonicity}.

        By \eqref{E.eq: Lbound2}, we have 
        \[
        q_{yy' \mid 1} = \pr(Y = y \mid Z = 1) + \pr(Y_0 = y' \mid Z = 1) - \{\pr(Y_0 \leq y \mid Z = 1) - \pr(Y_1 < y' \mid Z = 1)\}.
        \] 
        Therefore,  
        \begin{align} \label{E.eq: Area2}
        \sum_{k=y',k \neq y}^{J-1}\sum_{\ell=0,\ell \neq y'}^{y} q_{k\ell \mid 1} = \sum_{k=y',k \neq y}^{J-1}\sum_{\ell=0,\ell \neq y'}^{y} \pr(Y_1 = k, Y_0 = \ell \mid Z = 1) = 0. 
        \end{align}
        This implies  
        \[
        q_{k\ell \mid 1} = 0
        \] 
        for $k = y', \ldots, y-1, y+1, \ldots, J-1$ and $\ell = 0,\ldots, y'-1,y'+1,\ldots, y$. 
        Then, we have
        \begin{align*}
            q_{ky' \mid 1} = q_{k+ \mid 1}, \quad q_{y\ell \mid 1} = q_{+\ell \mid 1}, 
            \quad (k = y', \ldots, y-1; \ell = y'+1, \ldots, y). 
        \end{align*}

        After determining the above elements, $Q_1$ is shown in Figure 1. 
        All elements in Area \uppercase\expandafter{\romannumeral1} are zeros due to Assumption \ref{assume::monotonicity}. 
        The elements in Area \uppercase\expandafter{\romannumeral2} are determined by \eqref{E.eq: Area2}. 
        Area \uppercase\expandafter{\romannumeral3} and Area \uppercase\expandafter{\romannumeral4} contain all the undetermined elements in $Q_1$.

        \begin{figure}[t]
            \centering
            \includegraphics[width=0.9\linewidth]{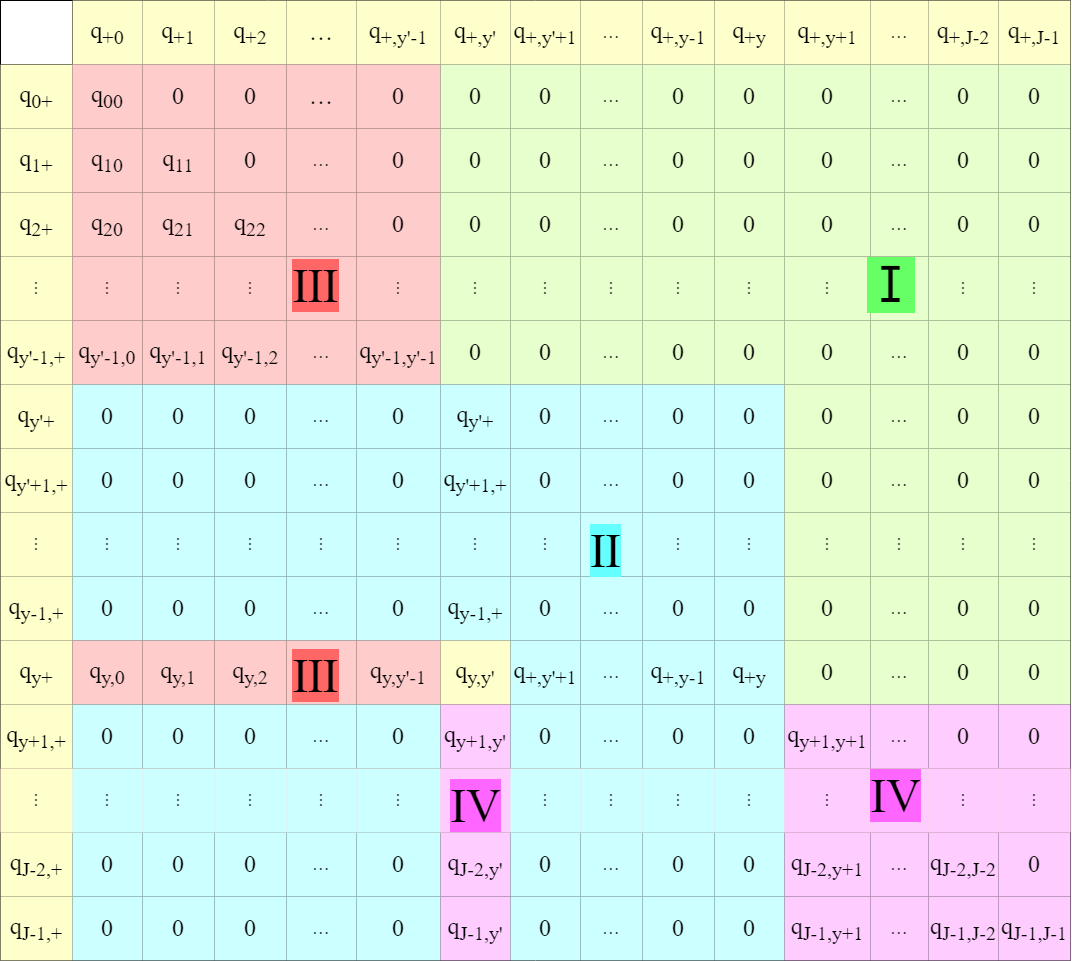}
            \caption{Visualization of the existence of the matrix $Q_1$ such that $\PN(Y_0=y', y)$ attains $\pr(Y = y \mid Z = 1) + \pr(Y_0 = y' \mid Z = 1) - \{\pr(Y_0 \leq y \mid Z = 1) - \pr(Y_1 < y' \mid Z = 1)\}$ under monotonicity}
            \label{fig:enter-label}
        \end{figure}

        The elements $q_{k\ell \mid 1}$ with $k,\ell \leq y'-1$ and $q_{y\ell \mid 1}$ with $\ell < y'-1$ in $Q_1$ form a submatrix: 
        \[ 
        Q_1^1=\left(\begin{array}{ccc}
        q_{00}^1 & \ldots  & q_{0,y'-1}^1  \\
        \vdots   & \ddots  & \vdots        \\
        q_{y',0}^1 &\ldots & q_{y',y'-1}^1 
        \end{array}\right) 
        =\left(\begin{array}{ccc}
        q_{00 \mid 1}     &\ldots  & q_{0,y'-1 \mid 1}  \\
           \vdots         &\ddots  & \vdots        \\
        q_{y'-1,0 \mid 1} &\ldots  & q_{y'-1, y'-1 \mid 1}  \\
        q_{y,0 \mid 1}    &\ldots  &q_{y,y'-1 \mid 1}
        \end{array}\right). 
        \]
        Then, the row sums and column sums of the submatrix $Q_1^1$ are 
        \begin{align*}
	   q^1_{1} &= \left(q^1_{0+}, \ldots, q^1_{y',+} \right)^\T
        = \left(q_{0+ \mid 1}, \ldots,q_{y'-1,+ \mid 1}, q_{y+ \mid 1} - q_{yy' \mid 1} - \sum_{\ell=y'+1}^{y} q_{+\ell \mid 1} \right)^\T, \\
	   q^1_{0} &= \left(q^1_{+0}, \ldots, q^1_{+,y'-1}\right)
        = \left(q_{+0 \mid 1}, \ldots, q_{+,y'-1 \mid 1}\right)^\T, 
        \end{align*}
        which satisfy $q^1_{k+} \geq 0$ for $k=0,\ldots,y'$ and $q^1_{+\ell} \geq 0$ for $\ell=0,\ldots,y'-1$. 
        In addition, 
        \[
        \sum_{k=0}^{y'} q^1_{k+} = \sum_{\ell=0}^{y'-1} q^1_{+\ell}= \sum_{k=0}^{y'}\sum_{\ell=0}^{y'-1} q^1_{k\ell}= \sum_{\ell=0}^{y'-1} q_{+\ell \mid 1}. 
        \]
        Apart from $q^1_{1}$ and $q^1_{0}$, there are no additional constraints on $Q_1^1$. 
        Therefore, by Lemma \ref{lemma: existence}, at least one matrix $Q^1_1$ satisfies the constraints of row sums $q^1_1$ and column sums $q^1_0$.

        The $q_{k\ell \mid 1}$'s with $k,\ell > y$ and $q_{ky' \mid 1}$ with $k > y$ in $Q_1$ form a submatrix: 
        \[ 
        Q_1^2=\left(\begin{array}{ccc}
        q_{00}^2      & \ldots  & q_{0,J-y-1}^2  \\
        \vdots        & \ddots  & \vdots        \\
        q_{J-y-2,0}^2 &\ldots & q_{J-y-2,J-y-1}^2 
        \end{array}\right) 
        =\left(\begin{array}{cccc}
        q_{y+1,y' \mid 1}  & q_{y+1,y+1 \mid 1}   &\ldots  & q_{y+1, J-1 \mid 1}  \\
           \vdots          & \vdots               &\ddots  & \vdots               \\
        q_{J-1,y' \mid 1}  & q_{J-1,y+1 \mid 1}   &\ldots  & q_{J-1, J-1 \mid 1} \\
        \end{array}\right), 
        \]
        where $q_{k\ell \mid 1} = 0$ for $k<\ell$. 
        Then, the row sums and column sums of the submatrix $Q_1^2$ are
        \begin{align*}
	   q^2_{1} &= \left(q^2_{0+}, \ldots, q^2_{J-y-1,+} \right)^\T
        = \left(q_{y+1,+ \mid 1}, \ldots,q_{J-1,+ \mid 1} \right)^\T, \\
         q^2_{0} &= \left(q^2_{+0}, \ldots, q^2_{+,J-y}\right)
        = \left(q_{+y' \mid 1} - q_{yy' \mid 1} - \sum_{k=y'}^{y-1} q_{k+ \mid 1}, q_{+,y+1 \mid 1}, \ldots, q_{+,J-1 \mid 1}\right)^\T, 
        \end{align*}
        which satisfy $q^2_{k+} \geq 0$ for $k=0,\ldots,J-y-1$ and $q^2_{+\ell} \geq 0$ for $\ell=0,\ldots,J-y$. 
        In addition, 
        \[
        \sum_{k=0}^{J-y-1} q^2_{k+} = \sum_{\ell=0}^{J-y} q^2_{+\ell}= \sum_{k=0}^{J-y-1}\sum_{\ell=0}^{J-y} q^2_{k\ell}= \sum_{k=y+1}^{J-1} q_{k+ \mid 1}. 
        \]
        Apart from $q^2_{1}$ and $q^2_{0}$, there are no additional constraints on $Q_1^2$. 
        By Lemma \ref{lemma: existence}, at least one matrix $Q^2_1$ satisfies the constraints of row sums $q^2_1$ and column sums $q^2_0$.  
    \end{itemize}

\noindent{\bf Step 2:} Proving the sharpness of the lower bound for $y'=y$.

    By \eqref{eq::D.lbound2}, the attainability of the lower bound implies that there exists a joint probability matrix $Q_1$ satisfying \eqref{E.eq: constraints PN} such that equality holds. 
    This means that if $\pr(Y = y \mid Z = 1) + \pr(Y_0 = y \mid Z = 1) - \{\pr(Y_0 \leq y \mid Z = 1) - \pr(Y_1 < y \mid Z = 1) \} \leq 0$, we have 
    \begin{equation}\label{E.eq: Lbound1'}
        \pr(Y_0=y, Y_1=y \mid Z=1) = 0.
    \end{equation}
    Otherwise, 
    \begin{align} \label{E.eq: Lbound2'}
        &\pr(Y_0=y, Y_1=y \mid Z=1)\\ \nonumber
        &= \pr(Y = y \mid Z = 1) + \pr(Y_0 = y \mid Z = 1) - \{\pr(Y_0 \leq y \mid Z = 1) - \pr(Y_1 < y \mid Z = 1)\}. 
    \end{align}

    Then, we prove the sharpness of the lower bound by showing the existence of a matrix $Q_1$ satisfying \eqref{E.eq: constraints PN} such that $\pr(Y_1=y,Y_0=y \mid Z=1)$ attains the bound in two parts.

    \begin{itemize}
        \item [1.] If $\pr(Y = y \mid Z = 1) + \pr(Y_0 = y \mid Z = 1) - \{\pr(Y_0 \leq y \mid Z = 1) - \pr(Y_1 < y \mid Z = 1) \} \leq 0$, $Q_1$ is subject to \eqref{E.eq: constraints PN} and \eqref{E.eq: Lbound1'}. 
        By \eqref{E.eq: Lbound1'}, we have 
        \begin{equation*}\label{E.eq: proof Lbound1'}
	\begin{aligned}
		  \pr(Y_1 = y,Y_0 = y \mid Z = 1) 
            = q_{yy \mid 1}
            = 0. 
	\end{aligned}
	\end{equation*}
        Then, we construct the elements of the $y$th row in $Q_1$ as
        \[
	       q_{y0 \mid 1} = \min\left\{q_{y+ \mid 1}, q_{+0 \mid 1}\right\}, \quad 
	       q_{y\ell \mid 1} = \min\left\{q_{y+ \mid 1} - \sum_{j=0}^{\ell-1} q_{yj \mid 1}, q_{+\ell \mid 1}\right\}, \quad
            (\ell = 1, \ldots, y-1).
        \]
        Let $Q_1^* = (q^*_{k\ell})_{k=0, \ldots, J-2; \ell=0, \ldots, J-1}$ denote the submatrix of $Q_1$ composed of all rows except the $y$th row. 
        Then, the row sums and column sums of the submatrix $Q_1^*$ are 
        \begin{align*}
	   q^*_{1} &= \left(q^*_{0+}, \ldots, q^*_{J-2,+ }\right)^\T
    	=\left(q_{0+ \mid 1}, \ldots,q_{y-1,+ \mid 1}, q_{y+1,+ \mid 1}, \ldots, q_{J-1,+ \mid 1}\right)^\T, \\
	   q^*_{0} &= \left(q^*_{+0}, \ldots, q^*_{+,J-1}\right)^\T
    	= \left(q_{+0 \mid 1} - q_{y0 \mid 1}, \ldots, q_{+,J-1 \mid 1} - q_{y,J-1 \mid 1}\right)^\T,
        \end{align*}
        which satisfy $q^*_{k+} \geq 0$ for $k=0,\ldots,J-2$ and $q^*_{+\ell} \geq 0$ for $\ell=0,\ldots,J-1$. 
         Because $\pr(Y = y \mid Z = 1) + \pr(Y_0 = y \mid Z = 1) - \{\pr(Y_0 \leq y \mid Z = 1) - \pr(Y_1 < y \mid Z = 1) \} \leq 0$, we have 
         \[
         \sum_{k=0}^{J-2} q^*_{k+} = \sum_{\ell=0}^{J-1} q^*_{+\ell}= \sum_{k=0}^{J-2}\sum_{\ell=0}^{J-1} q^*_{k\ell \mid 1}= 1 - q_{y+ \mid 1}. 
         \]
        Apart from $q^*_{1}$ and $q^*_{0}$, there are no additional constraints on $Q_1^*$. 
        By Lemma \ref{lemma: existence}, at least one matrix $Q^*_1$ satisfies the constraints of row sums $q^*_1$ and column sums $q^*_0$ which implies the existence of $Q_1$.

        \item [2.] If $\pr(Y = y \mid Z = 1) + \pr(Y_0 = y \mid Z = 1) - \{\pr(Y_0 \leq y \mid Z = 1) - \pr(Y_1 < y \mid Z = 1) \} > 0$, $Q_1$ is subject to \eqref{E.eq: constraints PN} and \eqref{E.eq: Lbound2'}. 
        In this part, we first construct the determinable elements in $Q_1$ based on \eqref{E.eq: constraints PN}, \eqref{E.eq: Lbound2'} and the monotonicity assumption. 
        Then, we demonstrate that the undetermined elements of $Q_1$ satisfy \eqref{E.eq: constraints PN}, \eqref{E.eq: Lbound2'}, and the monotonicity assumption.

        By \eqref{E.eq: Lbound2'}, we have 
        \[q_{yy \mid 1} = \pr(Y = y \mid Z = 1) + \pr(Y_0 = y \mid Z = 1) - \{\pr(Y_0 \leq y \mid Z = 1) - \pr(Y_1 < y \mid Z = 1)\}.\] 
        Therefore,  
        \begin{align}\label{E.eq: Area2'}
        \sum_{k=y+1}^{J-1}\sum_{\ell=0}^{y-1} q_{k\ell \mid 1} = \sum_{k=y+1}^{J-1}\sum_{\ell=0}^{y-1} \pr(Y_1 = k, Y_0 = \ell \mid Z = 1) = 0. 
        \end{align}
        This implies  
        \[q_{k\ell \mid 1} = 0\] 
        for $k = y+1, \ldots, J-1$ and $\ell = 0,\ldots, y-1$. 
        Under Assumption \ref{assume::monotonicity}, we have
        \begin{align*}
            q_{k\ell \mid 1} = 0, 
            \quad (k < \ell). 
        \end{align*}

        After determining the above elements, $Q_1$ is shown in Figure 2. 
        All elements in Area \uppercase\expandafter{\romannumeral1} are zeros due to Assumption \ref{assume::monotonicity}. 
        The elements in Area \uppercase\expandafter{\romannumeral2} are determined by \eqref{E.eq: Area2'}. 
        Area \uppercase\expandafter{\romannumeral3} and Area \uppercase\expandafter{\romannumeral4} contain all the undetermined elements in $Q_1$.
        \begin{figure}{t}
            \centering
            \includegraphics[width=0.7 \linewidth]{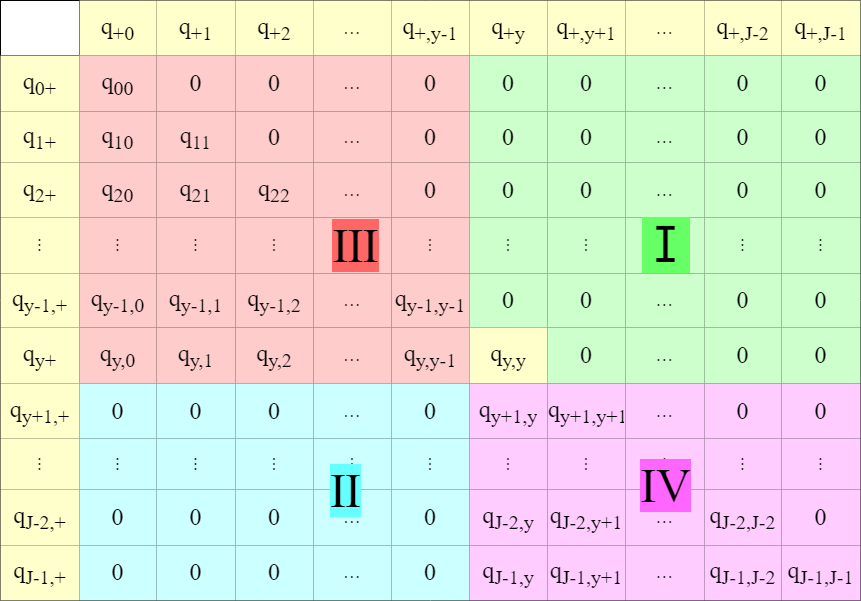}
            \caption{Visualization of the existence of the matrix $Q_1$ such that $\PN(Y_0=y, y)$ attains $\Pr(Y = y \mid Z = 1) + \Pr(Y_0 = y \mid Z = 1) - \{\Pr(Y_0 \leq y \mid Z = 1) - \Pr(Y_1 < y \mid Z = 1)\}$ under monotonicity}
            \label{fig2}
        \end{figure}

        The $q_{k\ell \mid 1}$'s with $k \leq y$ and $\ell < y$ in $Q_1$ form a submatrix: 
        \[ 
        Q_1^1=\left(\begin{array}{ccc}
        q_{00}^1 & \ldots  & q_{0,y-1}^1  \\
        \vdots   & \ddots  & \vdots        \\
        q_{y,0}^1 &\ldots & q_{y,y-1}^1 
        \end{array}\right) 
        =\left(\begin{array}{ccc}
        q_{00 \mid 1}     &\ldots  & q_{0,y-1 \mid 1}  \\
           \vdots         &\ddots  & \vdots        \\
        q_{y,0 \mid 1} &\ldots  & q_{y, y-1 \mid 1} 
        \end{array}\right). 
        \]
        Then, the row sums and column sums of the submatrix $Q_1^1$ are 
        \begin{align*}
	   q^1_{1} &= \left(q^1_{0+}, \ldots, q^1_{y,+} \right)^\T
        = \left(q_{0+ \mid 1}, \ldots,q_{y-1,+ \mid 1}, q_{y+ \mid 1} - q_{yy \mid 1} \right)^\T, \\
	   q^1_{0} &= \left(q^1_{+0}, \ldots, q^1_{+,y-1}\right)
        = \left(q_{+0 \mid 1}, \ldots, q_{+,y-1 \mid 1}\right)^\T, 
        \end{align*}
        which satisfy $q^1_{k+} \geq 0$ for $k=0,\ldots,y$ and $q^1_{+\ell} \geq 0$ for $\ell=0,\ldots,y-1$. 
        In addition, 
        \[
        \sum_{k=0}^{y} q^1_{k+} = \sum_{\ell=0}^{y-1} q^1_{+\ell}= \sum_{k=0}^{y}\sum_{\ell=0}^{y-1} q^1_{k\ell}= \sum_{\ell=0}^{y-1} q_{+\ell \mid 1}. 
        \]
        Apart from $q^1_{1}$ and $q^1_{0}$, there are no additional constraints on $Q_1^1$. 
        Therefore, by Lemma \ref{lemma: existence}, at least one matrix $Q^1_1$ satisfies the constraints of row sums $q^1_1$ and column sums $q^1_0$.

        The $q_{k\ell \mid 1}$'s with $k > y$ and $\ell \geq y$ in $Q_1$ form a submatrix: 
        \[ 
        Q_1^2=\left(\begin{array}{ccc}
        q_{00}^2      & \ldots  & q_{0,J-y-1}^2  \\
        \vdots        & \ddots  & \vdots        \\
        q_{J-y-2,0}^2 &\ldots & q_{J-y-2,J-y-1}^2 
        \end{array}\right) 
        =\left(\begin{array}{cccc}
        q_{y+1,y \mid 1}  &\ldots  & q_{y+1, J-1 \mid 1}  \\
           \vdots         &\ddots  & \vdots               \\
        q_{J-1,y \mid 1}  &\ldots  & q_{J-1, J-1 \mid 1} \\
        \end{array}\right), 
        \]
        where $q_{k\ell \mid 1} = 0$ for $k<\ell$. 
        Then, the row sums and column sums of the submatrix $Q_1^2$ are
        \begin{align*}
	   q^2_{1} &= \left(q^2_{0+}, \ldots, q^2_{J-y-1,+} \right)^\T
        = \left(q_{y+1,+ \mid 1}, \ldots,q_{J-1,+ \mid 1} \right)^\T, \\
	   q^2_{0} &= \left(q^2_{+0}, \ldots, q^2_{+,J-y}\right)
        = \left(q_{+y \mid 1} - q_{yy \mid 1}, q_{+,y+1 \mid 1} \ldots, q_{+,J-1 \mid 1}\right)^\T, 
        \end{align*}
        which satisfy $q^2_{k+} \geq 0$ for $k=0,\ldots,J-y-1$ and $q^2_{+\ell} \geq 0$ for $\ell=0,\ldots,J-y$. 
        In addition, 
        \[
        \sum_{k=0}^{J-y-1} q^2_{k+} = \sum_{\ell=0}^{J-y} q^2_{+\ell}= \sum_{k=0}^{J-y-1}\sum_{\ell=0}^{J-y} q^2_{k\ell}= \sum_{k=y+1}^{J-1} q_{k+ \mid 1}. 
        \]
        Apart from $q^2_{1}$ and $q^2_{0}$, there are no additional constraints on $Q_1^2$. 
        By Lemma \ref{lemma: existence}, at least one matrix $Q^2_1$ satisfies the constraints of row sums $q^2_1$ and column sums $q^2_0$.  
    \end{itemize}
\end{proof}

\subsection{Proving the sharpness of the upper bound in Theorem \ref{theorem::bounds-mono}} \label{subsection: bounds-mono-upper}

In this subsection, we prove the sharpness of the upper bounds in two steps.

\begin{proof}
\noindent{\bf Step 1:} Proving the sharpness of the upper bound for $y'<y$.

    By \eqref{eq::D.ubound1}, we prove the sharpness of the upper bound by showing the existence of a matrix $Q_1$ satisfying \eqref{E.eq: constraints PN} such that $\PN(\omega_0,y)$ attains the bound in three parts.

    \begin{itemize}
        \item [1.] If $\min_{y' \leq h < y} \left\{\pr(Y_1 = y \mid Z = 1), \pr(Y_0 = y' \mid Z = 1), \pr(Y_0 \leq h \mid Z = 1) - \pr(Y_1 \leq h \mid Z = 1) \right\} = \pr(Y = y \mid Z = 1)$, we have
        \begin{equation} \label{E.eq: Ubound1}
            \pr(Y_1=y, Y_0=y' \mid Z=1) = \pr(Y = y \mid Z = 1). 
        \end{equation} 
        We next prove the existence of a matrix $Q_1$ satisfying \eqref{E.eq: constraints PN} and \eqref{E.eq: Ubound1}. 
        By \eqref{E.eq: Ubound1}, we have 
        \begin{equation*}
            \pr(Y_1=y, Y_0 \neq y' \mid Z=1) = \sum_{\ell=0, \ell \neq y'}^{J-1} q_{y\ell \mid 1} = 0.
        \end{equation*}
        This implies 
         \[
        q_{yy' \mid 1} = q_{y+ \mid 1}, \quad q_{y\ell \mid 1} = 0, \quad (\ell = 0, \ldots,y'-1, y'+1, \ldots, J-1). 
        \]
        Let $Q_1^*$ denote the submatrix of $Q_1$ composed of all rows except the $y$th row. 
        Then, the row sums and column sums of the submatrix $Q_1^*$ are 
        \begin{align*}
	   q^*_{1} &= \left(q^*_{0+}, \ldots, q^*_{J-2,+ }\right)^\T
    	=\left(q_{0+ \mid 1}, \ldots,q_{y-1,+ \mid 1}, q_{y+1,+ \mid 1}, \ldots, q_{J-1,+ \mid 1}\right)^\T, \\
        q^*_{0} &= \left(q^*_{+0}, \ldots, q^*_{+,J-1}\right)^\T
    	= \left(q_{+0 \mid 1} - q_{y0 \mid 1}, \ldots, q_{+,J-1 \mid 1} - q_{y,J-1 \mid 1}\right)^\T, 
        \end{align*}
        which satisfy $q^*_{k+} \geq 0$ for $k=0,\ldots,J-2$ and $q^*_{+\ell} \geq 0$ for $\ell=0,\ldots,J-1$. 
        In addition, we have 
        \[
	\sum_{k=0}^{J-2} q^*_{k+}
	= \sum_{\ell=0}^{J-1} q^*_{+\ell}
	= \sum_{k=0}^{J-2}\sum_{\ell=0}^{J-1} q^*_{k\ell}
	= 1- q_{y+ \mid 1}.
        \]
        Apart from $q^*_{1}$ and $q^*_{0}$, there are no additional constraints on $Q_1^*$. 
        By Lemma \ref{lemma: existence}, at least one matrix $Q^*_1$ satisfies the constraints of row sums $q^*_1$ and column sums $q^*_0$ which implies the existence of $Q_1$.

        \item [2.] If $\min_{y' \leq h < y} \{\pr(Y_1 = y \mid Z = 1), \pr(Y_0 = y' \mid Z = 1), \pr(Y_0 \leq h \mid Z = 1) - \pr(Y_1 \leq h \mid Z = 1) \} = \pr(Y_0 = y' \mid Z = 1)$, we have
        \begin{equation} \label{E.eq: Ubound2}
            \pr(Y_1=y,Y_0=y' \mid Z=1) = \pr(Y_0=y' \mid Z=1). 
        \end{equation}
        We next prove the existence of a matrix $Q_1$ satisfying \eqref{E.eq: constraints PN} and \eqref{E.eq: Ubound2}. 
        By \eqref{E.eq: Ubound2}, we have 
        \[
        \pr(Y_1 \neq y, Y_0=y' \mid Z=1) = \sum_{k=0,k\neq y}^{J-1} q_{ky' \mid 1} = 0, 
        \]
        which implies 
        \[
        q_{yy' \mid 1} = q_{y+ \mid 1}, \quad q_{ky' \mid 1} = 0, \quad (k = 0, \ldots, y-1, y+1, \ldots, J-1). 
        \]
        Let $Q_1^* = (q^*_{k\ell})_{k=0, \ldots, J-1; \ell=0, \ldots, J-2}$ denote the submatrix of $Q_1$ composed of all columns except the $y$th column. 
        Then, the row sums and column sums of the submatrix $Q_1^*$ are 
        \begin{align*}
	   q^*_{1} &= \left(q^*_{0+}, \ldots, q^*_{J-1,+ }\right)^\T
    	=\left(q_{0+ \mid 1} - q_{0y' \mid 1}, \ldots, q_{J-1,+ \mid 1} - q_{J-1,y' \mid 1}\right)^\T, \\
	   q^*_{0} &= \left(q^*_{+0}, \ldots, q^*_{+,J-2}\right)^\T
    	= \left(q_{+0 \mid 1}, \ldots, q_{+,y'-1 \mid 1}, q_{+,y'+1 \mid 1}, \ldots, q_{+,J-1 \mid 1}\right)^\T, 
        \end{align*}
        which satisfy $q^*_{k+} \geq 0$ for $k=0,\ldots,J-1$ and $q^*_{+\ell} \geq 0$ for $\ell=0,\ldots,J-2$. 
        In addition, 
        \[
	\sum_{k=0}^{J-1} q^*_{k+}
	= \sum_{\ell=0}^{J-2} q^*_{+\ell}
	= \sum_{k=0}^{J-1}\sum_{\ell=0}^{J-2} q^*_{k\ell}
	= 1- q_{+y' \mid 1}.
        \]
        Apart from $q^*_{1}$ and $q^*_{0}$, there are no additional constraints on $Q_1^*$. 
        By Lemma \ref{lemma: existence}, at least one matrix $Q^*_1$ satisfies the constraints of row sums $q^*_1$ and column sums $q^*_0$ which implies the existence of $Q_1$.

         \item [3.] 
         If $\min_{y' \leq h < y} \{\pr(Y_1 = y \mid Z = 1), \pr(Y_0 = y' \mid Z = 1), \pr(Y_0 \leq h \mid Z = 1) - \pr(Y_1 \leq h \mid Z = 1)\} = \pr(Y_0 \leq h^* \mid Z = 1) - \pr(Y_1 \leq h^* \mid Z = 1)$, we have
        \begin{equation}\label{E.eq: Ubound3}
            \pr(Y_1=y,Y_0=y' \mid Z=1) = \pr(Y_0 \leq h^* \mid Z = 1) - \pr(Y_1 \leq h^* \mid Z = 1), 
        \end{equation}
        where $h^* \in \{y',\ldots,y-1\}$ and 
        $
        \min_{y' \leq h < y} \{\pr(Y_0 \leq h \mid Z = 1) - \pr(Y_1 \leq h \mid Z = 1)\} = \pr(Y_0 \leq h^* \mid Z = 1) - \pr(Y_1 \leq h^* \mid Z = 1). 
        $

        We next prove the existence of a matrix $Q_1$ satisfying \eqref{E.eq: constraints PN} and \eqref{E.eq: Ubound3}. 
        We represent $Q_1$ into block matrix based on $h^*$: 
        \[ 
        Q=\left(\begin{array}{cc}
        Q_1^1 &\quad  Q_1^2  \\
        Q_1^3 &\quad  Q_1^4 
        \end{array}\right), 
        \]
        where 
        \begin{align*}
       & Q_1^1=\left(\begin{array}{ccc}
        q_{00 \mid 1}  &\quad \ldots &\quad  q_{0h^* \mid 1}  \\
            \vdots     &\quad \ddots &\quad      \vdots      \\
        q_{h^*0 \mid 1} &\quad \ldots &\quad  q_{h^*h^* \mid 1}
        \end{array}\right), \  \qquad
        Q_1^2=\left(\begin{array}{ccc}
        q_{0,h^*+1 \mid 1}  &\quad \ldots &\quad  q_{0,J-1 \mid 1}  \\
            \vdots        &\quad \ddots &\quad      \vdots      \\
        q_{h^*,h^*+1 \mid 1}  &\quad \ldots &\quad  q_{h^*,J-1 \mid 1}
        \end{array}\right), \\
       & Q_1^3=\left(\begin{array}{ccc}
        q_{h^*+1,0 \mid 1}  &\quad \ldots &\quad  q_{h^*+1,h \mid 1}  \\
            \vdots        &\quad \ddots &\quad      \vdots      \\
        q_{J-1,0 \mid 1}  &\quad \ldots &\quad  q_{J-1,h^* \mid 1}
        \end{array}\right), \quad 
        Q_1^4=\left(\begin{array}{ccc}
        q_{h^*+1,h^*+1 \mid 1} &\quad \ldots &\quad  q_{h^*+1,J-1 \mid 1}  \\
            \vdots         &\quad \ddots &\quad      \vdots        \\
        q_{J-1,h^*+1 \mid 1} &\quad \ldots &\quad  q_{J-1,J-1 \mid 1}
        \end{array}\right).  
        \end{align*}

        We can determine all elements in $Q_1^2$ and $Q_1^3$, based on constraints \eqref{E.eq: constraints PN}, \eqref{E.eq: Ubound3} and the monotonicity assumption.

        Under Assumption \ref{assume::monotonicity}, $Q_1^2$ is a zero matrix, and $Q_1^1$ and $Q_1^4$ are lower triangular matrixes.

        In $Q_1^3$, by $\pr(Y_1 = y, Y_0 = y' \mid Z = 1) = \pr(Y_0 \leq h^* \mid Z = 1) - \pr(Y_1 \leq h^* \mid Z = 1)$, we have 
        \[
        q_{yy' \mid 1} = \sum_{k=h^*}^{J-1}\sum_{\ell=0}^{h^*} q_{k\ell \mid 1}.
        \]
        This equation implies 
        \[
        q_{yy' \mid 1} = \sum_{\ell \leq h^*} q_{+\ell \mid 1} - \sum_{k \leq h^*} q_{k+ \mid 1}, \quad q_{k\ell \mid 1} = 0, 
        \]
        for $\ell=0,\ldots, y'-1, y'+1, \ldots h^*$ and $k = h^*+1, \ldots y-1, y+1, \ldots, J-1$. 
        Therefore, all elements in $Q_1^3$ are determined.

        Then, we prove the existence of $Q_1$ by showing the existence of $Q_1^1$ and $Q_1^4$ satisfy the constraints. 
        For $Q_1^1$, because $y' \leq h^* < y$, we have the row sums and column sums are 
        \begin{align*}
	   q^1_{1} &= \left(q^1_{0+}, \ldots, q^1_{h^*+}\right)^\T
            = \left(q_{0+ \mid 1}, \ldots, q_{h^*+ \mid 1}\right)^\T, \\
	   q^1_{0} &= \left(q^1_{+0}, \ldots, q^1_{+h^*} \right)^\T = \left(q_{+0 \mid 1} - q_{y0 \mid 1}, \ldots, q_{+h^* \mid 1} - q_{yh^* \mid 1} \right)^\T, 
        \end{align*}
        which satisfy $q^1_{k+} \geq 0$ and $q^1_{+\ell} \geq 0$ for $k,\ell=0,\ldots,h^*$. 
        In addition, 
        \[
        \sum_{k=0}^{h^*} q^1_{k+} = \sum_{\ell=0}^{h^*} q^1_{+\ell} = \sum_{k=0}^{h^*} \sum_{\ell=0}^{h^*} q^1_{k\ell} = \sum_{k=0}^{h^*} q_{k+ \mid 1}. 
        \]
        Apart from $q^1_{1}$ and $q^1_{0}$, there are no additional constraints on $Q^1_1$. 
        By Lemma \ref{lemma: existence}, at least one matrix $Q^1_1$ satisfies the constraints of row sums $q^1_1$ and column sums $q^1_0$.

        Similarly, we have the row sums and column sums for $Q_1^4$ are 
        \begin{align*}
        q_1^4 &= \left(q^4_{0+}, \ldots, q^4_{J-h^*-2,+}\right)
        = \left(q_{h^*+1,+ \mid 1} - q_{h^*+1,y' \mid 1}, \ldots, q_{J-1,+ \mid 1} - q_{J-1,y' \mid 1}\right), \\ 
        q_0^4 &= \left(q^4_{+0}, \ldots, q^4_{+,J-h^*-2}\right)
        = \left(q_{+,h^*+1 \mid 1}, \ldots, q_{+,J-1 \mid 1}\right), 
        \end{align*}
        which satisfy $q^4_{k+} \geq 0$ and $q^4_{+\ell} \geq 0,$ for $k,\ell=h^*+1,\ldots,J-1$.
        In addition, 
        \[
        \sum_{k=0}^{J-h^*-2} q^4_{k+} = \sum_{\ell=0}^{J-h^*-2} q^4_{+\ell} = \sum_{k=0}^{J-h^*-2} \sum_{\ell=0}^{J-h^*-2} q^4_{k\ell} = \sum_{k=h^*+1}^{J-1} q_{k+ \mid 1}. 
        \]
        Apart from $q^4_{1}$ and $q^4_{0}$, there are no additional constraints on $Q^4_1$. 
        By Lemma \ref{lemma: existence}, at least one matrix $Q^4_1$ satisfies the constraints of row sums $q^4_1$ and column sums $q^4_0$. 
    
        Therefore, there exists a matrix $Q_1$ such that $\pr(Y_1 = y, Y_0 = y' \mid Z = 1) = \pr(Y_0 \leq h^* \mid Z = 1) - \pr(Y_1 \leq h^* \mid Z = 1)$ under \eqref{E.eq: constraints PN}. 
    \end{itemize}

\noindent{\bf Step 2:} Proving the sharpness of the upper bound for $y'=y$.

By \eqref{eq::D.ubound2}, we prove the sharpness of the upper bound by showing the existence of a matrix $Q_1$ satisfying \eqref{E.eq: constraints PN} such that $\pr(Y_1=y,Y_0=y \mid Z=1)$ attains the bound in two parts.

    \begin{itemize}
        \item [1.] If $\min \left\{\pr(Y_1 = y \mid Z = 1), \pr(Y_0 = y \mid Z = 1) \right\} = \pr(Y = y \mid Z = 1)$, we have
        \begin{equation} \label{E.eq: Ubound1'}
            \pr(Y_1=y, Y_0=y \mid Z=1) = \pr(Y = y \mid Z = 1). 
        \end{equation} 
        We next prove the existence of a matrix $Q_1$ satisfying \eqref{E.eq: constraints PN} and \eqref{E.eq: Ubound1'}. 
        By \eqref{E.eq: Ubound1'}, we have 
        \begin{equation*}
            \pr(Y_1=y, Y_0 \neq y \mid Z=1) = \sum_{\ell=0}^{y-1} q_{y\ell \mid 1} = 0.
        \end{equation*}
        This implies 
         \[
        q_{yy \mid 1} = q_{y+ \mid 1}, \quad q_{y\ell \mid 1} = 0, \quad (\ell = 0, \ldots, y-1). 
        \]
        Let $Q_1^*$ denote the submatrix of $Q_1$ composed of all rows except the $y$th row. 
        Then, the row sums and column sums of the submatrix $Q_1^*$ are 
        \begin{align*}
	   q^*_{1} &= \left(q^*_{0+}, \ldots, q^*_{J-2,+ }\right)^\T
    	=\left(q_{0+ \mid 1}, \ldots,q_{y-1,+ \mid 1}, q_{y+1,+ \mid 1}, \ldots, q_{J-1,+ \mid 1}\right)^\T, \\
	   q^*_{0} &= \left(q^*_{+0}, \ldots, q^*_{+,J-1}\right)^\T
    	= \left(q_{+0 \mid 1} - q_{y0 \mid 1}, \ldots, q_{+,J-1 \mid 1} - q_{y,J-1 \mid 1}\right)^\T, 
        \end{align*}
        which satisfy $q^*_{k+} \geq 0$ for $k=0,\ldots,J-2$ and $q^*_{+\ell} \geq 0$ for $\ell=0,\ldots,J-1$. 
        In addition, we have 
        \[
	\sum_{k=0}^{J-2} q^*_{k+}
	= \sum_{\ell=0}^{J-1} q^*_{+\ell}
	= \sum_{k=0}^{J-2}\sum_{\ell=0}^{J-1} q^*_{k\ell}
	= 1- q_{y+ \mid 1}.
        \]
        Apart from $q^*_{1}$ and $q^*_{0}$, there are no additional constraints on $Q_1^*$. 
        By Lemma \ref{lemma: existence}, at least one matrix $Q^*_1$ satisfies the constraints of row sums $q^*_1$ and column sums $q^*_0$ which implies the existence of $Q_1$.

        \item [2.] If $\min \{\pr(Y_1 = y \mid Z = 1), \pr(Y_0 = y \mid Z = 1)\} = \pr(Y_0 = y \mid Z = 1)$, we have
        \begin{equation} \label{E.eq: Ubound2'}
            \pr(Y_1=y,Y_0=y \mid Z=1) = \pr(Y_0=y \mid Z=1). 
        \end{equation}
        We next prove the existence of a matrix $Q_1$ satisfying \eqref{E.eq: constraints PN} and \eqref{E.eq: Ubound2'}. 
        By \eqref{E.eq: Ubound2'}, we have 
        \[
        \pr(Y_1 \neq y, Y_0=y \mid Z=1) = \sum_{k=y+1}^{J-1} q_{ky \mid 1} = 0, 
        \]
        which implies 
        \[
        q_{yy \mid 1} = q_{y+ \mid 1}, \quad q_{ky \mid 1} = 0, \quad (k = y+1, \ldots, J-1). 
        \]
        Let $Q_1^* = (q^*_{k\ell})_{k=0, \ldots, J-1; \ell=0, \ldots, J-2}$ denote the submatrix of $Q_1$ composed of all columns except the $y$th column. 
        Then, the row sums and column sums of the submatrix $Q_1^*$ are 
        \begin{align*}
	   q^*_{1} &= \left(q^*_{0+}, \ldots, q^*_{J-1,+ }\right)^\T
    	=\left(q_{0+ \mid 1} - q_{0y \mid 1}, \ldots, q_{J-1,+ \mid 1} - q_{J-1,y \mid 1}\right)^\T, \\
	   q^*_{0} &= \left(q^*_{+0}, \ldots, q^*_{+,J-2}\right)^\T
    	= \left(q_{+0 \mid 1}, \ldots, q_{+,y-1 \mid 1}, q_{+,y+1 \mid 1}, \ldots, q_{+,J-1 \mid 1}\right)^\T, 
        \end{align*}
        which satisfy $q^*_{k+} \geq 0$ for $k=0,\ldots,J-1$ and $q^*_{+\ell} \geq 0$ for $\ell=0,\ldots,J-2$. 
        In addition, 
        \[
	\sum_{k=0}^{J-1} q^*_{k+}
	= \sum_{\ell=0}^{J-2} q^*_{+\ell}
	= \sum_{k=0}^{J-1}\sum_{\ell=0}^{J-2} q^*_{k\ell}
	= 1- q_{+y \mid 1}.
        \]
        Apart from $q^*_{1}$ and $q^*_{0}$, there are no additional constraints on $Q_1^*$. 
        By Lemma \ref{lemma: existence}, at least one matrix $Q^*_1$ satisfies the constraints of row sums $q^*_1$ and column sums $q^*_0$ which implies the existence of $Q_1$. 
    \end{itemize} 
\end{proof}

\subsection{Numerical solutions via linear programming} \label{subsection: linear programming}
We can obtain the sharp bounds on general $\PN(\omega_0, y)$ under monotonicity using linear programming similar to \cite{tianpearl2000}. 
By \eqref{C.eq: defPN}, we have 
\[
\PN(\omega_0,y)
=\frac{\pr(Y_1=y, \omega_0 \mid Z=1)}{\pr(Y=y \mid Z=1)}
=\frac{\sum_{\ell=0}^{J-1} c_\ell q_{y\ell \mid 1}}{\pr(Y=y \mid Z=1)} 
=\frac{\Tilde{c}^{\T}q}{\pr(Y=y \mid Z=1)}, 
\]
where $q=(q_{00\mid 1},\ldots, q_{0,J-1 \mid 1}, q_{10 \mid 1}, \ldots, q_{1,J-1 \mid 1}, \ldots,q_{J-1,0 \mid 1}, \ldots, q_{J-1,J-1 \mid 1})^{\T}$ and $\Tilde{c}$ is a binary column vector of length $J^2$, determined by $\omega_0$. 
This implies that we only need to derive the sharp bounds on $\pr(Y_1=y, \omega_0 \mid Z=1) = \Tilde{c}^{\T}q$ under constraints in \eqref{eq: constraints PN} 
and constraints provided by the monotonicity assumption. 
Then, the sharp bounds on $\Tilde{c}^{\T}q$ can be obtained by solving the following linear programming: 
\begin{equation}\label{C.eq: prime}
\begin{aligned}
	\max _q \quad &\Tilde{c}^{\T}q \\
     \text{ s.t. } Aq &=b \\
 	              q &\geq 0;
\end{aligned}
\end{equation}
where 
\begin{itemize}
    \item  $A=(A_1, A_2)^\T$: $A_1^\T$ is a numeric matrix of dimensions $(2J-1)\times J^2$, consisting of 0 and 1, and uniquely determined by the constraints in \eqref{eq: constraints PN}, while $A_2^\T$ is a numeric matrix of dimensions $J(J-1)/2 \times J^2$, consisting of 0 and 1, and uniquely determined by the monotonicity assumption;
    \item $b=(b_1,b_2)^\T$: $b_1^\T = \left(\pr(Y_1=0 \mid Z=1), \ldots, \pr(Y_1=J-2 \mid Z=1), \pr(Y_0=0 \mid Z=1), \ldots, \pr(Y_0=J-\right. \\ \left. 2 \mid Z=1), 1\right)$ is determined by the constraints in \eqref{eq: constraints PN}, while $b_2^\T$ is a column vector of dimension $J(J-1)/2$ with all elements being zero determined by the monotonicity.
\end{itemize}

Solving \eqref{C.eq: prime} involves symbolic computation. For simplicity, we can consider its dual problem: 
\begin{equation}\label{C.eq: dual}
\begin{aligned}
	\min _q \quad &b^{\T} \lambda \\
	\text{s.t. } A^{\T} \lambda &\geq \Tilde{c}. 
\end{aligned}
\end{equation}
In \eqref{C.eq: dual}, the minimum of $b^{\T} \lambda$ can be obtained by enumerating all the vertices of a convex polygon $A^{\T} \lambda \geq \Tilde{c}$. 
By the properties of linear programming, if the dual problem has an optimal solution $\lambda^{*}$, then the primal has an optimal solution $P^{*}$, and $\Tilde{c}^{\T} P^{*}$ = $b^{\T} \lambda^{*}$. Therefore, we can obtain the sharp upper bound on $\pr(Y_1=y, \omega_0 \mid Z=1)$ under constraints in \eqref{eq: constraints PN}. 
Similarly, we can obtain the sharp lower bound on $\pr(Y_1=y,\omega_0 \mid Z=1)$.

\section{Proof of Proposition \ref{proposition::falsify}} \label{section: proposition}

Recall the definition of $\delta_{k \mid 1}$ in \eqref{eq: delta}. 
By Lemma \ref{lemma::identification-Q}, by requiring $q_{k\ell \mid 1}$'s to satisfy the Fr\'echet bounds in Lemma \ref{lemma: frechet}, we have
    \begin{align*}
    \max \left(\begin{array}{c}
			0 \\
			\pr(Y_1 = k \mid Z = 1) + \pr(Y_0 = k - 1 \mid Z = 1) - 1
		\end{array}\right)
        \leq &q_{k,k-1 \mid 1} = \delta_{k \mid 1} \\
        \leq &\min \left(\begin{array}{c}
			\pr(Y_1 = k \mid Z = 1) \\
			\pr(Y_0 = k-1 \mid Z = 1)
		\end{array}\right), \\
     \max \left(\begin{array}{c}
			0 \\
			\pr(Y_1 = k \mid Z = 1) + \pr(Y_0 = k \mid Z = 1) - 1
		\end{array}\right)
        \leq &q_{k,k \mid 1} = \pr(Y_1 = k \mid Z = 1) - \delta_{k \mid 1}  \\
        \leq &\min \left(\begin{array}{c}
			\pr(Y_1 = k \mid Z = 1) \\
			\pr(Y_0 = k \mid Z = 1)
		\end{array}\right),
    \end{align*}
    where $k = 1, \ldots, J-1$. Therefore, 
    \begin{align*}
		\max \left(\begin{array}{c}
			0 \\
			\pr(Y_1 = k \mid Z = 1) + \pr(Y_0 = k - 1 \mid Z = 1) - 1
		\end{array}\right)
        \leq &\delta_{k \mid 1} \leq
        \min \left(\begin{array}{c}
			\pr(Y_1 = k \mid Z = 1) \\
			\pr(Y_0 = k-1 \mid Z = 1)
		\end{array}\right), 
    \end{align*}
     for $z = 0, 1$ and $k = 1, \ldots, J-1$.

\section{Results on the probability of causation} \label{section: discussion}

\subsection{Definition of the probability of causation} \label{subsection::PC-definition}

\cite{dawidfaigmanfienberg2014} proposed the definition of the probability of causation for binary outcomes as $\PC = \pr\left(Y_0 = 0 \mid Y_{1} = 1\right)$. 
We propose a general definition for the probability of causation for ordinal outcomes. 
Let $\omega_0$ denote an event only depending on $Y_0$, we define the probability of causation as 
\[
    \PC(\omega_0, y) = \pr(\omega_0 \mid Y_1 = y). 
\]
We will derive results on $\PC(\omega_0,y)$  in parallel with those on $\PN(\omega_0,y)$ in the main paper. 
The probability of causation only involves the joint distribution of the potential outcomes and is independent of the treatment assignment mechanism. Therefore, we focus on randomized experiments.

\begin{assumption}[Randomization] \label{assume::randomization}
$Z$ is independent of $(Y_1,Y_0)$.  
\end{assumption}

Under Assumption \ref{assume::randomization}, we can identify $\pr(Y_1)$ and $\pr(Y_0)$. 
Since the definition of $\PC(\omega_0, y)$ involves the joint distribution $\pr(Y_1, Y_0)$, $\PC(\omega_0,y)$ is still not identifiable without additional assumptions.

\subsection{Main results on the probability of causation} \label{subsection::PC-results}

Let $P = (p_{k\ell})_{k,\ell= 0,\ldots, J-1}$ where $p_{k\ell} = \pr(Y_1=k, Y_0=\ell)$. The $p_{k\ell}$'s satisfy the following constraints: 
\begin{equation} \label{D.eq: constraints PC}
\begin{aligned}
    &\sum_{\ell=0}^{J-1} p_{k\ell}=\pr(Y_1=k),\quad 
     \sum_{k=0}^{J-1} p_{k\ell}=\pr(Y_0=\ell), \\
    &\sum_{k=0}^{J-1}\sum_{\ell=0}^{J-1} p_{k\ell}=1,\quad  
      p_{k\ell} \geq 0, \quad \left(k, \ell=0, \ldots, J-1 \right).
\end{aligned}
\end{equation}

Similar to the identification of $q_{k\ell \mid 1}$'s in the main paper, we have the following lemma.

\begin{lemma}  \label{lemma: identify P} 
    Under Assumption \ref{assume::incremental}, $P$ is identified by 
    \begin{align*}
	p_{00} &= \pr(Y_1 = 0),\\
         p_{k,k-1} &= \sum_{j=0}^{k-1} \{ \pr(Y_0 = j) - \pr(Y_1 = j)\},\\
	p_{kk} &= \sum_{j=0}^{k} \pr(Y_1 = j) - \sum_{j=0}^{k-1} \pr(Y_0 = j),\\
         p_{k'\ell'} &= 0.
    \end{align*}
     for $k = 1, \ldots, J-1$ and $k' < \ell'$ or $k' > \ell'+1$.   
\end{lemma}

$\PC(\omega_0, y)$ is a linear combination of $p_{k\ell}$'s as 
\begin{align} \label{F.eq: defPC}
    \PC(\omega_0, y) = \frac{\sum_{\ell=0}^{J-1} c_\ell p_{y\ell}}{\pr(Y_1 = y)}, 
\end{align}
where $c_{\ell}$'s are binary constants and uniquely determined by $\omega_0$. 
Therefore, we can identify $\PC(\omega_0, y)$ if we can identify $P$. 
Lemma \ref{lemma: identify P} implies the following theorem, with the definition of 
\begin{equation}
\label{eq::xi-k-appendix}
\xi_k = \sum_{j=0}^{k-1} \{\pr(Y_0 = j) - \pr(Y_1 = j)\}.
\end{equation}

\begin{theorem} \label{th: identify PC}
    Under Assumptions \ref{assume::randomization} and \ref{assume::incremental}, $\PC(\omega_0,y)$ is identified by
    \[
	\PC(\omega_0,y) = c_{y} + \frac{(c_{y-1} - c_y)\xi_y }{\pr(Y_1 = y)}, 
    \]
    where $c_{y-1}$ and $c_y$ are uniquely determined by $\omega_0$ in \eqref{F.eq: defPC}.
\end{theorem}

Assumption \ref{assume::incremental} is key to the identification result in Theorem \ref{th: identify PC}. There are testable implications based on Lemma \ref{lemma: identify P} because the $p_{k\ell}$'s must satisfy the Fr\'echet bounds \citep{ruschendorf1991frechet}. We present the result in Proposition \ref{proposition::falsify-pc} below.

\begin{proposition} \label{proposition::falsify-pc}
 Under Assumption \ref{assume::marginal}, Assumption \ref{assume::incremental} implies
    \begin{align*}
  	&\max \left(\begin{array}{c}
			0,
			\pr(Y = k \mid Z=1) + \pr(Y = k - 1 \mid Z=0) - 1
		\end{array}\right) \\
        &\leq \xi_k \leq 
        \min \left(\begin{array}{c}
			\pr(Y = k \mid Z=1),
			\pr(Y = k-1 \mid Z=0)
		\end{array}\right),
    \end{align*}
     for $k = 1, \ldots, J-1$ and $z = 0, 1$, where $\xi_{k}$ is defined in \eqref{eq::xi-k-appendix}.
\end{proposition}

If the inequality in Proposition \ref{proposition::falsify-pc} fails, then Assumption \ref{assume::incremental} is falsified. However, Assumption \ref{assume::incremental} cannot be validated by data. Even if the inequality in Proposition \ref{proposition::falsify-pc} holds, Assumption \ref{assume::incremental} can still fail.

When Assumption \ref{assume::incremental} does not hold, $\PC(\omega_0, y)$ is generally not identifiable. Without Assumption \ref{assume::incremental}, we can derive the sharp bounds on $\PC(\omega_0, y)$ under Assumption \ref{assume::marginal}.

\begin{theorem} \label{th: bounds PC}
    Under Assumption \ref{assume::randomization}, the sharp bounds on $\PC(\omega_0, y)$ are
    \begin{align*}
	\max \left(\begin{array}{c}
		0 \\
		\frac{\pr(Y_1 = y) - \sum_{\ell=0}^{J-1} (1 - c_\ell) \pr(Y_0 = \ell)}{\pr(Y_1 = y)}
	\end{array}\right)
        \leq \PC(\omega_0, y) \leq
         \min \left(\begin{array}{c}
		1 \\
		\frac{\sum_{\ell=0}^{J-1} c_{\ell}\pr(Y_0 = \ell) }{\pr(Y_1 = y)}\\
	\end{array}\right),
    \end{align*}
    where $c_\ell$'s are determined by $\omega_0$ in \eqref{F.eq: defPC}. 
\end{theorem}

If we are willing to maintain the monotonicity in Assumption \ref{assume::monotonicity}, then we can derive narrower sharp bounds on $\pr(Y_0 \neq y \mid Y_1 = y)$ and $\pr(Y_0 = y' \mid Y_1 = y)$.

\begin{theorem}
\label{theorem::pc-bounds-mono}
Under Assumptions \ref{assume::randomization} and \ref{assume::monotonicity}, the sharp bounds on $\pr(Y_0 \neq y \mid Y_1 = y)$ are    
$$
\max \left(\begin{array}{c}
			0,
			\frac{\pr(Y_1 = y) - \pr(Y_0 = y)}{\pr(Y_1 = y)}
		\end{array}\right)
        \leq \pr(Y_0 \neq y \mid Y_1 = y) \leq
        \min \left(\begin{array}{c}
			1,
			\frac{\xi_{y}}{\pr(Y_1 = y)}
		\end{array}\right),
$$
and the sharp bounds on $\pr(Y_0 = y' \mid Y_1 = y)$ are 
\begin{align*}
		&\max \left(\begin{array}{c}
			0,
			\frac{\pr(Y_1 = y) + \sum_{k=0}^{y'-1} \pr(Y_1 = k) - \sum_{\ell = 0}^y \pr(Y_0 = \ell) + \pr(Y_0 = y')}{\pr(Y_1 = y)}
		\end{array}\right) \\
        &\leq \pr(Y_0 = y' \mid Y_1 = y) \leq
        \min_{y' < k \leq y} \left(\begin{array}{c}
			1,
            \frac{\pr(Y_0 = y')}{\pr(Y_1 = y)},
			\frac{\xi_{k}}{\pr(Y_1 = y)}
		\end{array}\right),
  \end{align*}
  recalling the definition of $\xi_k$ in \eqref{eq::xi-k-appendix}. 
\end{theorem}

The proofs of the above results on $\PC(\omega_0, y)$ are in parallel with those on $\PN(\omega_0,y)$. We only need to modify the conditional probabilities given $Z=1$ as the unconditional probabilities. To avoid repetition, we omit the details.

\bibliographystyle{apalike}
\bibliography{ordinalY}

\end{document}